\documentclass[journal]{IEEEtran}

\usepackage{graphicx}
\usepackage{color}
\usepackage{amsfonts}
\usepackage{float}
\usepackage{amsmath}
\usepackage{amsthm}
\usepackage{epstopdf}
\usepackage{changepage}
\usepackage{latexsym}
\usepackage{tikz}
\usepackage{algorithm,algorithmic}
\usepackage{amssymb}
\usepackage{multicol}
\usepackage{cases}
\usepackage{url}
\usepackage{xcolor}
\usetikzlibrary{arrows}
\usetikzlibrary{automata}
\usetikzlibrary{er}
\usetikzlibrary{positioning}
\usetikzlibrary{shadows}

\newcommand{\tabincell}[2]{\begin{tabular}{@{}#1@{}}#2\end{tabular}}

\floatstyle{plain}
\newfloat{myalgo}{tbhp}{mya}
{\begin{myalgo}[#1]
\centering
\begin{minipage}{#2}
\begin{algorithm}[H]}%
{\end{algorithm}
\end{minipage}
\end{myalgo}}

\DeclareMathOperator*{\mycup}{\bigcup}

\begin{document}

\title{Optimal PMU Placement for Power System Dynamic State Estimation by Using Empirical Observability Gramian}

\author{Junjian~Qi,~\IEEEmembership{Member,~IEEE,}
        Kai~Sun,~\IEEEmembership{Senior Member,~IEEE,}
        and Wei~Kang,~\IEEEmembership{Fellow,~IEEE}
        \thanks{This work was supported in part by the University of Tennessee, Knoxville, the CURENT Engineering Research Center, and NRL. Paper no. TPWRS-00394-2014.

        J.~Qi and K. Sun are with the Department of Electrical Engineering and Computer Science, University of Tennessee, Knoxville, TN 37996 USA (e-mails: junjian.qi.2012@ieee.org; kaisun@utk.edu).

        W. Kang is with the Department of Applied Mathematics, Naval Postgraduate School, Monterey, CA 93943 USA (e-mail: wkang@nps.edu).
}
}

\maketitle
\markboth{preprint of doi 10.1109/TPWRS.2014.2356797, IEEE Transactions on Power Systems.}{stuff}\maketitle

\begin{abstract}
In this paper the empirical observability Gramian calculated around the operating region of a power system is used to quantify the degree of observability of the system states under specific phasor measurement unit (PMU) placement. An optimal PMU placement method for power system dynamic state estimation is further formulated as an optimization problem which maximizes the determinant of the empirical observability Gramian and is efficiently solved by the NOMAD solver, which implements the Mesh Adaptive Direct Search (MADS) algorithm. The implementation, validation, and also the robustness to load fluctuations and contingencies of the proposed method are carefully discussed. The proposed method is tested on WSCC 3-machine 9-bus system and NPCC 48-machine 140-bus system by performing dynamic state estimation with square-root unscented Kalman filter. The simulation results show that the determined optimal PMU placements by the proposed method can guarantee good observability of the system states, which further leads to smaller estimation errors and larger number of convergent states for dynamic state estimation compared with random PMU placements. Under optimal PMU placements an obvious observability transition can be observed. The proposed method is also validated to be very robust to both load fluctuations and contingencies.
\end{abstract}


\section{Introduction} \label{intro}

\IEEEPARstart{P}{ower} system state estimation is an important application of
the energy management system (EMS). The most widely studied static state estimation [\ref{se1}]--[\ref{se6}]
cannot capture the dynamics of power systems very well due to its dependency on slow update rates of the Supervisory Control and Data Acquisition (SCADA) systems.
Accurate dynamic states of the system obtained from real-time dynamic state estimation
facilitated by high-level phasor measurement unit (PMU) deployment has thus become essential.
Because of the high global positioning system (GPS) synchronization accuracy, PMUs can provide highly synchronized direct measurements of voltage phasors and current phasors, thus playing a critical role in achieving real-time wide-area monitoring, protection, and control.
If the dynamic states of the system can be tracked in real time, wide-area monitoring and control schemes that require synchronized system-wide phasor data could be further implemented against major stability problems [\ref{sun1}], [\ref{sun2}].

One important question for dynamic state estimation is the observability of nonlinear power systems.
A well-selected subset of variables can contain sufficient information about the rest of the variables,
thus allowing us to reconstruct the system's complete internal states and to make the system observable.
Therefore, it is important to install the PMUs at appropriately chosen buses
so as to maximize the observability of the system states and further efficiently achieve the desired functionalities.

The observability of a system can be checked by the observability rank condition [\ref{kailath}]--[\ref{diop1}].
But this only offers a yes or no answer and is limited to small systems due to heavy computation burden. Alternatively, the empirical observability Gramian [\ref{lall}]--[\ref{Gramianref}] provides a computable tool for empirical analysis of the state-output behaviour, which has been used in various applications of control system observabilities [\ref{Krener}], [\ref{kang1}]. Recently it has also been applied to power systems [\ref{sun}].
Furthermore, methodologies of optimal sensor placement based on observability Gramian have been developed for weather prediction [\ref{kang}] and chemical engineering [\ref{singh}]--[\ref{det}].

Most research on PMU placement is for static state estimation and is mainly based on the topological observability criterion, which only specifies that the power system states should be uniquely estimated with the minimum number of PMU measurements but neglects important parameters such as transmission line admittances by only focusing on the binary connectivity graph [\ref{ob1}], [\ref{info}].
Under this framework, many approaches have been proposed, such as mixed integer programming [\ref{ob2}], [\ref{ob3}], binary search [\ref{ob4}], metaheuristics [\ref{ob5}], [\ref{ob6}], particle swarm optimization [\ref{ob7}], and eigenvalue-eigenvector based approaches [\ref{ob8}], [\ref{ob9}]. Different from the conventional approaches, an information-theoretic criterion, the mutual information between PMU measurements and power system states, is proposed to generate highly informative PMU configurations [\ref{info}].

By contrast, not much research has been undertaken on PMU placement for power system dynamic state estimation.
Numerical PMU configuration algorithms have been proposed to maximize the overall sensor response
while minimizing the correlation among sensor outputs [\ref{kamwa}].
The system observability for a PMU placement scheme and the corresponding uncertainty are evaluated via the steady-state error covariance obtained from discrete algebraic Riccati equation [\ref{zhang1}], [\ref{zhang2}].
A PMU placement strategy for dynamic state estimation has also been proposed to
ensure a satisfactory state tracking performance [\ref{huang}].

In this paper, the empirical observability Gramian [\ref{lall}]--[\ref{Gramianref}] is applied to
quantify the degree of observability of the system states
and formulate the optimal PMU placement for power system
dynamic state estimation as an optimization problem
that maximizes the determinant of the empirical observability Gramian.

The remainder of this paper is organized as follows.
Section \ref{s_ob} introduces the fundamentals of observability and the definition of
empirical observability Gramian.
Section \ref{pmu placement} discusses the formulation of optimal PMU placement, the generator and measurement model, and the implementation, validation, and robustness of the proposed method.
In Section \ref{case} the proposed optimal PMU placement method is tested and validated on WSCC 3-machine 9-bus system and NPCC 48-machine 140-bus system. Finally the conclusion is drawn in Section \ref{conclusion}.

\section{Fundamentals of Observability and Empirical Observability Gramian} \label{s_ob}

A system is observable if the system's complete internal state can be reconstructed from its outputs.
For a linear time-invariant system
\begin{subnumcases} {\label{linear}}
\dot{\boldsymbol{x}}=A\boldsymbol{x}+B\boldsymbol{u} \\
\boldsymbol{y}=C\boldsymbol{x}+D\boldsymbol{u}
\end{subnumcases}
where $\boldsymbol{x} \in \mathbb{R}^n$ is the state vector and $\boldsymbol{y}\in \mathbb{R}^p$ is the output vector,
it is observable if the observability matrix
\begin{displaymath}
\left[ \begin{array}{c}
C \\
CA \\
CA^2 \\
\vdots \\
CA^{n-1}
\end{array} \right]
\end{displaymath}
or the observability Gramian [\ref{kailath}]
\begin{equation}
\boldsymbol{W}_{o,\textrm{linear}}=\int_0^\infty e^{A^Tt}C^TCe^{At}dt
\end{equation}
has full rank.

For a nonlinear system
\begin{subnumcases} {\label{n1}}
\dot{\boldsymbol{x}}=\boldsymbol{f}(\boldsymbol{x},\boldsymbol{u}) \\
\boldsymbol{y}=\boldsymbol{h}(\boldsymbol{x},\boldsymbol{u})
\end{subnumcases}
where $\boldsymbol{f}(\cdot)$ and $\boldsymbol{h}(\cdot)$ are the state transition and output functions, $\boldsymbol{x} \in \mathbb{R}^n$ is the state vector, $\boldsymbol{u} \in \mathbb{R}^v$ is the input vector, and $\boldsymbol{y}\in \mathbb{R}^p$ is the output vector,
it is locally observable at a state $\boldsymbol{x}_0$ if the nonlinear observability matrix obtained by using Lie derivative
has full rank at $\boldsymbol{x}=\boldsymbol{x}_0$ [\ref{diop}], [\ref{diop1}].

The rank test method is easy and straightforward for linear systems and it
can tell if the system is observable under a specific sensor set.
However, for nonlinear systems this can be very complicated even for small systems.

One possibility is to linearize the nonlinear system. But the nonlinear dynamics of the system will be lost.
Alternatively, empirical observability Gramian [\ref{lall}], [\ref{lall1}] provides a computable tool for
empirical analysis of the state-output behaviour of a nonlinear system.
It is also proven that the empirical observability Gramian of a stable linear system
described by (\ref{linear}) is equal to the usual observability Gramian [\ref{lall1}].
Singh and Hahn [\ref{singh}]--[\ref{det}] show that it can be used for observability analysis of
nonlinear systems over an operating region and can be readily computed for systems of considerable size.

The following sets are defined for empirical observability Gramian:
\begin{align}
T^n&=\{T_1,\cdots,T_r;\;\;\;T_l \in \mathbb{R}^{n\times n},\;T_l^T T_l=I_n,\;l=1,\cdots,r\} \nonumber \\
M&=\{c_1,\cdots,c_s;\;\;\;\;c_m \in \mathbb{R},\;c_m>0,\;m=1,\cdots,s\} \nonumber \\
E^n&=\{e_1,\cdots,e_n;\;\;\;\textrm{standard unit vectors in}\;\mathbb{R}^n\} \nonumber
\end{align}
where $T^n$ defines the initial state perturbation directions, $r$ is the number of matrices for perturbation directions, $I_n$ is an identity matrix with dimension $n$, $M$ defines the perturbation sizes and $s$ is the number of different perturbation sizes for each direction;
and $E^n$ defines the state to be perturbed and $n$ is the number of states of the system.

For the nonlinear system described by (\ref{n1}), the empirical observability Gramian can be defined as
\begin{equation} \label{gd0}
\boldsymbol{W}=\sum_{l=1}^{r}\sum_{m=1}^{s}\frac{1}{rsc_m^2}\int_0^\infty T_l \Psi^{lm}(t)T_l^T dt
\end{equation}
where $\Psi^{lm}(t)\in \mathbb{R}^{n\times n}$ is given by $\Psi_{ij}^{lm}(t)=(y^{ilm}(t)-y^{ilm,0})^T(y^{jlm}(t)-y^{jlm,0})$, $y^{ilm}(t)$ is the output of the nonlinear system corresponding to the initial condition $\boldsymbol{x}(0)=c_mT_le_i+\boldsymbol{x}_0$, and $y^{ilm,0}$ refers to the output measurement corresponding to the unperturbed initial state $\boldsymbol{x}_0$, which is usually chosen as the steady state under typical power flow conditions but can also be chosen as other operating points.

In practical implementation, (\ref{gd0}) can be rewritten as its discrete form [\ref{Gramianref}]
\begin{equation} \label{gd}
\boldsymbol{W}=\sum_{l=1}^{r}\sum_{m=1}^{s}\frac{1}{rsc_m^2}\sum_{k=0}^K T_l \Psi^{lm}_k T_l^T \Delta t_k
\end{equation}
where $\Psi^{lm}_k \in \mathbb{R}^{n\times n}$ is given by ${\Psi^{lm}_k}_{ij}=(y^{ilm}_k-y^{ilm,0})^T(y^{jlm}_k-y^{jlm,0})$, $y^{ilm}_k$ is the output at time step $k$,  $K$ is the number of points chosen for the approximation of the integral in (\ref{gd0}), and $\Delta t_k$ is the time interval between two points.

For multiple outputs $\boldsymbol{y}\in \mathbb{R}^p$, ${\Psi_k^{lm}}_{ij}$ in (\ref{gd}) is
\renewcommand{\arraystretch}{1.6}
\[ {\Psi_k^{lm}}_{ij}=\left[ \begin{array}{c}
y_{1,k}^{ilm}-y_1^{ilm,0} \\
y_{2,k}^{ilm}-y_2^{ilm,0} \\
\vdots \\
y_{p,k}^{ilm}-y_p^{ilm,0} \end{array} \right]^T
\left[ \begin{array}{c}
y_{1,k}^{jlm}-y_1^{jlm,0} \\
y_{2,k}^{jlm}-y_2^{jlm,0} \\
\vdots \\
y_{p,k}^{jlm}-y_p^{jlm,0} \end{array} \right]
\]
\[\quad\; =\sum\limits_{o=1}^p (y_{o,k}^{ilm}-y_o^{ilm,0})^T (y_{o,k}^{jlm}-y_o^{jlm,0}) \]
where $y_{o,k}^{ilm}$ and $y_o^{ilm,0}$ are respectively the $o$th output of the nonlinear system
corresponding to the perturbed and unperturbed initial state.

Then the matrix $\Psi^{lm}_k$ for $p$ outputs can be given as
\begin{equation} \label{psi}
\Psi^{lm}_k=\sum\limits_{o=1}^p\Psi_{o,k}^{lm}
\end{equation}
where $\Psi_{o,k}^{lm}$ is calculated for output $o$.

By substituting (\ref{psi}) into (\ref{gd}), we can get
\begin{equation}
\boldsymbol{W}=\sum\limits_{o=1}^p \boldsymbol{W}^o
\end{equation}
where $\boldsymbol{W}^o$ is the empirical Gramian
for ouput $o$ and can be calculated by substituting $\Psi_{o,k}^{lm}$ into (\ref{gd}).

Therefore, the empirical observability Gramian for a system with $p$ outputs is the summation of the empirical Gramians computed for each of the $p$ outputs individually [\ref{singh1}].

Different from analysis based on linearization, the empirical observability Gramian is defined using the original system model. It reflects the observability of the full nonlinear dynamics in the given domain, whereas the observability based on linearization only works locally in a neighborhood of an operating point.

\section{Optimal PMU Placement for Power System Dynamic State Estimation} \label{pmu placement}

In this section the optimal PMU placement for dynamic state estimation is formulated based on the empirical observability Gramian. The generator and measurement model and the implementation, validation, and robustness of the proposed method are also discussed.

\subsection{Formulation of Optimal PMU Placement} \label{formulate}

The degree of observability can be quantified by making use of a variety of different measures of the empirical observability Gramian, such as the smallest eigenvalue [\ref{Krener}]--[\ref{kang}], [\ref{muller}], the trace [\ref{singh1}], the determinant [\ref{det}], [\ref{muller}], or the condition number [\ref{Krener}].
Although all of them are based on Gramian matrices, different measures reflect various aspects of observability.
The smallest eigenvalue defines the worst scenario of observability.
It measures the largest possible error among all dimensions of unknown noises.
The trace of Gramian matrices measures the total gain from state variation to sensor output.
It cannot detect the unobservability of individual directions in noise space.
Observability based on the condition number emphasizes the numerical stability in state estimation.
The determinant of Gramian matrices measures the overall observability in all directions in noise space.

Although the trace of Gramian also tends to measure the overall observability,
it cannot tell the existence of a zero eigenvalue. Thus an unobservable system
may still have a large trace. Compared with the method based on the smallest eigenvalue,
the determinant is a smooth function, which is a desirable property in numerical computations.
However, it is advised that the smallest eigenvalue be verified to be at an acceptable level
when using the determinant of Gramian as the measure of observability.
This is to avoid the situation in which a Gramian has an almost zero eigvenvalue
that makes the system practically unobservable. In fact, we numerically verify
the smallest eigenvalue in all the examples.
Based on all these considerations, in this paper we choose the objective as the maximization of the determinant of
the observability Gramian under different PMU placements.

To better understand the PMU placement method based on the maximization of the determinant of the empirical observability Gramian, consider the state trajectory, $\boldsymbol{x}(t,\boldsymbol{x}_0)$ of the system with an initial state $\boldsymbol{x}_0$. The corresponding output function of sensor measurement is $\boldsymbol{y}(t,\boldsymbol{x}_0)$. Their relationship can be treated as a mapping from $\boldsymbol{x}_0$ to $\boldsymbol{y}(t,\boldsymbol{x}_0)$. The process of state estimation is to find the inverse mapping from $\boldsymbol{y}(t,\boldsymbol{x}_0)$ to $\boldsymbol{x}(t,\boldsymbol{x}_0)$. It is well known that an inverse problem can be solved accurately, or is well posed, if the image of a mapping is sensitive to the variation of the input variable. Thus it is desirable to have a large gain from $\boldsymbol{x}_0$ to $\boldsymbol{y}(t,\boldsymbol{x}_0)$. Since the empirical observability Gramian is approximately the gain from $\boldsymbol{x}_0$ to $\boldsymbol{y}(t,\boldsymbol{x}_0)$ [\ref{kang1}], selecting sensor locations with a larger determinant of the Gramian improves the overall observability.

Based on this choice, the optimal PMU placement problem can be formulated as
\begin{align} \label{opt}
&\max\limits_{\boldsymbol{z}} \det \, \boldsymbol{W}(\boldsymbol{z}) \nonumber \\
\textrm{s.t.}\; & \;\;\; \sum_{i=1}^g  z_i = \bar{g} \\
& \;\;\; z_i \in \{0,1\}, \;\; i=1,\cdots,g \nonumber
\end{align}
where $\boldsymbol{z}$ is the vector of binary control variables which determines where to place the PMUs, $\boldsymbol{W}$ is the corresponding empirical observability Gramian, $g$ is the number of generators in the system, and $\bar{g}$ is the number of PMUs to be placed.

\subsection{Simplified Generator and Measurement Model} \label{sim model}

A classical second-order generator model and simplified measurement model
are used for demonstration.
All the nodes except for the internal generator nodes are eliminated and the admittance matrix of
the reduced network $\boldsymbol{\overline{Y}}$ is obtained. The equations of motion of the generators are given by
\begin{subnumcases}{\label{gen model1}}
\dot{\delta_i}=\omega_i-\omega_0 \\
\dot{\omega}_i=\frac{\omega_0}{2H_i}(T_{mi}-T_{ei})
\end{subnumcases}
where $i=1,\cdots,g$, $\delta_i$, $\omega_i$, $T_{mi}$, $T_{ei}$, $\omega_0$, and $H_i$ are respectively the rotor angle, rotor speed in rad/s, mechanical torque, electric air-gap torque, rated value of angular frequency, and inertia of the $i$th generator; the electric air-gap torque $T_{ei}$ can be written as
\begin{align}
T_{ei}&\cong P_{ei}=E_i^2G_{ii} \nonumber \\
&+\sum_{\substack{j=1 \\j\ne i}}^{g}E_iE_j\big(G_{ij}\textrm{cos}(\delta_i-\delta_j)+B_{ij}\textrm{sin}(\delta_i-\delta_j)\big)
\end{align}
where $E_i$ and $E_j$ are voltage magnitudes of the $i$th and $j$th internal generator bus;
$G_{ii}$ and $G_{ij}$ are real elements of $\boldsymbol{\overline{Y}}$ and $B_{ij}$ is the imaginary element of $\boldsymbol{\overline{Y}}$.

For simplicity it is assumed that for generators where PMUs are installed the rotor angle and rotor speed can be directly measured. $T_{mi}$ and $E_i$ are used as inputs and are actually kept constant in the simulation. The dynamic model (\ref{gen model1}) can thus be rewritten in a general state space form in (\ref{n1}) and the state vector $\boldsymbol{x}$, input vector $\boldsymbol{u}$, and output vector $\boldsymbol{y}$
can be written as
\begin{subequations}
\begin{align}
\boldsymbol{x} &= [\boldsymbol{\delta}^T \quad \boldsymbol{\omega}^T]^T \\
\boldsymbol{u} &= [\boldsymbol{T_m}^T \quad \boldsymbol{E}^T]^T \\
\boldsymbol{y} &= [{\boldsymbol{\delta}^{\mathcal{G}_P}}^T {\boldsymbol{\omega}^{\mathcal{G}_P}}^T]^T
\end{align}
\end{subequations}
where $\mathcal{G}_{P}$ is the set of generators where PMUs are installed and $\boldsymbol{\delta}^{\mathcal{G}_P}$ and $\boldsymbol{\omega}^{\mathcal{G}_P}$ are the rotor angle and rotor speed of the generators that belong to $\mathcal{G}_P$.

The simplified generator and measurement model in this section is denoted by $\textrm{\textbf{M}}_1$, for which the number of states $n=2\,g$.

\subsection{Realistic Generator and Measurement Model}

In this section a more realistic generator and measurement model is presented by mainly following [\ref{zhou}].
There are some differences between our model and that in [\ref{zhou}].
Firstly, [\ref{zhou}] only consider the single-machine infinite-bus system but our model can be used for multi-machine systems. Secondly, [\ref{zhou}] only consider fourth-order transient generator model but we allow both the fourth-order transient model and the second-order classical model.
Thirdly, in [\ref{zhou}] only the terminal voltage phasor is used as output and the terminal current phasor is used as input but we consider both the terminal voltage phasor and terminal current phasor as outputs.
Lastly, [\ref{zhou}] does not require the admittance matrix for estimation while our model requires this knowledge.

Let $\mathcal{G}_4$ and $\mathcal{G}_2$ respectively denote the set of generators with fourth-order model and second-order model. The numbers of generators with fourth-order model or second-order model, which are also the cardinality of the set $\mathcal{G}_4$ and $\mathcal{G}_2$, are respectively $g_4$ and $g_2$.
Thus the number of states $n=4\,g_4+2\,g_2$.
For generator $i\in \mathcal{G}_4$, the fast sub-transient dynamics and saturation effects are ignored and the generator model is described by the fourth-order differential equations in local $d$-$q$ reference frame:
\begin{subnumcases}{\label{gen model}}
\dot{\delta_i}=\omega_i-\omega_0 \\
\dot{\omega}_i=\frac{\omega_0}{2H_i}(T_{mi}-T_{ei}-\frac{K_{Di}}{\omega_0}(\omega_i-\omega_0)) \\
\dot{e}'_{qi}=\frac{1}{T'_{d0i}}(E_{fdi}-e'_{qi}-(x_{di}-x'_{di})i_{di}) \\
\dot{e}'_{di}=\frac{1}{T'_{q0i}}(-e'_{di}+(x_{qi}-x'_{qi})i_{qi})
\end{subnumcases}
where $i$ is the generator serial number, $\delta_i$ is the rotor angle,
$\omega_i$ is the rotor speed in $rad/s$, and $e'_{qi}$ and $e'_{di}$ are the transient voltage along $q$ and $d$ axes; $i_{qi}$ and $i_{di}$ are stator currents at $q$ and $d$ axes;
$T_{mi}$ is the mechanical torque, $T_{ei}$ is the electric air-gap torque, and $E_{fdi}$ is the internal field voltage; $\omega_0$ is the rated value of angular frequency, $H_i$ is the inertia constant, and $K_{Di}$ is the damping factor; $T'_{q0i}$ and $T'_{d0i}$ are the open-circuit time constants for $q$ and $d$ axes; $x_{qi}$ and $x_{di}$ are the synchronous reactance and $x'_{qi}$ and $x'_{di}$ are the transient reactance respectively at the $q$ and $d$ axes.

For generator $i\in \mathcal{G}_2$, the generator model is only described by the first two equations of (\ref{gen model}) and $e'_{qi}$ and $e'_{di}$ are kept unchanged. Similar to Section \ref{sim model}, the set of generators where PMUs are installed is denoted by $\mathcal{G}_{P}$. For generator $i \in \mathcal{G}_{P}$, $T_{mi}$, $E_{fdi}$, the terminal voltage phaosr $E_{ti}=e_{Ri}+je_{Ii}$, and the terminal current phasor $I_{ti}=i_{Ri}+ji_{Ii}$ can be measured, among which $T_{mi}$ and $E_{fdi}$ are used as inputs and $E_{ti}$ and $I_{ti}$ are the outputs.

The dynamic model (\ref{gen model}) can be rewritten in a general state space form in (\ref{n1}) and
the state vector $\boldsymbol{x}$, input vector $\boldsymbol{u}$, and output vector $\boldsymbol{y}$
can be written as
\begin{subequations}
\begin{align}
\boldsymbol{x} &= [\boldsymbol{\delta}^T \quad \boldsymbol{\omega}^T \quad \boldsymbol{e'_{q}}^T \quad \boldsymbol{e'_{d}}^T]^T \\
\boldsymbol{u} &= [\boldsymbol{T_m}^T \quad \boldsymbol{E_{fd}}^T]^T \\
\boldsymbol{y} &= [\boldsymbol{e_R}^T \quad \boldsymbol{e_I}^T \quad \boldsymbol{i_R}^T \quad \boldsymbol{i_I}^T]^T.
\end{align}
\end{subequations}

The $i_{qi}$, $i_{di}$, and $T_{ei}$ in (\ref{gen model}) can be written as functions of $\boldsymbol{x}$ and $\boldsymbol{u}$:
\begin{subequations} \label{temp}
\begin{align}
\it \Psi_{Ri}&=e'_{di}\sin\delta_i+e'_{qi}\cos\delta_i \\
\it \Psi_{Ii}&=e'_{qi}\sin\delta_i-e'_{di}\cos\delta_i \\
I_{ti}&=\overline{Y}_i(\it \Psi_{R}+j\Psi_{I}) \\
i_{Ri}&= \operatorname{Re}(I_{ti}) \\
i_{Ii}&= \operatorname{Im}(I_{ti}) \\
i_{qi}&=i_{Ii}\sin\delta_i+i_{Ri}\cos\delta_i \\
i_{di}&=i_{Ri}\sin\delta_i-i_{Ii}\cos\delta_i \\
e_{qi}&=e'_{qi}-x'_{di}i_{di} \\
e_{di}&=e'_{di}+x'_{qi}i_{qi} \\
T_{ei}&\cong P_{ei}=e_{qi}i_{qi}+e_{di}i_{di}.
\end{align}
\end{subequations}
where $\it \Psi_i=\Psi_{Ri}+j\Psi_{Ii}$ is the voltage source, $\it \Psi_R$ and $\it \Psi_I$ are column vectors of all generators' $\it \Psi_{Ri}$ and $\it \Psi_{Ii}$, $e_{qi}$ and $e_{di}$ are the terminal voltage at $q$ and $d$ axes, and $\overline{Y}_i$ is the $i$th row of the admittance matrix of the reduced network $\boldsymbol{\overline{Y}}$.

In \eqref{temp} the outputs $i_R$ and $i_I$ are written as functions of $\boldsymbol{x}$ and $\boldsymbol{u}$.
Similarly, the outputs $e_{Ri}$ and $e_{Ii}$ can also be written as function of $\boldsymbol{x}$ and $\boldsymbol{u}$:
\begin{subequations}
\begin{align}
e_{Ri}&= e_{di}\sin\delta_i+e_{qi}\cos\delta_i \\
e_{Ii}&=e_{qi}\sin\delta_i-e_{di}\cos\delta_i.
\end{align}
\end{subequations}

Compared with the simplified generator measurement model $\textrm{\textbf{M}}_1$ in Section \ref{sim model},
the generator and measurement model considered here is more realistic, which is denoted by $\textrm{\textbf{M}}_2$.

Similar to [\ref{zhou}], the continuous models in (\ref{gen model1}) and (\ref{gen model}) can be discretized into their discrete form as
\begin{subnumcases} {\label{n2}}
\boldsymbol{x}_k = \boldsymbol{f}_d(\boldsymbol{x}_{k-1},\boldsymbol{u}_{k-1}) \\
\boldsymbol{y}_k = \boldsymbol{h}(\boldsymbol{x}_k,\boldsymbol{u}_k)
\end{subnumcases}
where $k$ denotes the time at $k\Delta_t$ and the state transition functions $\boldsymbol{f}_d$ can be obtained by the modified Euler method [\ref{kunder}] as
\begin{align}
\tilde{\boldsymbol{x}}_k &= \boldsymbol{x}_{k-1}+\boldsymbol{f}(\boldsymbol{x}_{k-1},\boldsymbol{u}_{k-1})\Delta t  \\
\tilde{\boldsymbol{f}} &= \frac{\boldsymbol{f}(\tilde{\boldsymbol{x}}_k,\boldsymbol{u}_k) + \boldsymbol{f}(\boldsymbol{x}_{k-1},\boldsymbol{u}_{k-1})}{2}\\
\boldsymbol{x}_k &= \boldsymbol{x}_{k-1}+\tilde{\boldsymbol{f}}\Delta t.
\end{align}

The model in (\ref{n2}) is used to obtain the system response and the outputs for the empirical observability Gramian calculation in (\ref{gd}).
It is also used to perform the dynamic state estimation with the squart-root unscented Kalman filter (SR-UKF) [\ref{sr_ukf}], which can further be used to validate the optimal PMU placement method proposed in this paper.

\subsection{Implementation} \label{impl}

In Section \ref{s_ob} it has been shown that the empirical observability Gramian for a system with $p$ outputs
is the summation of the empirical Gramians computed for each of the $p$ outputs individually.
This can be easily generalized to the case for placing $\bar{g}$ PMUs in power system dynamic state estimation.
The empirical observability Gramian calculated from placing $\bar{g}$ PMUs individually
adds to be the identical Gramian calculated from placing the $\bar{g}$ PMUs simultaneously, which can be shown as
\begin{equation}
\hspace*{-3.6cm} \boldsymbol{W}=\sum\limits_{o=1}^{p}\boldsymbol{W}^o =\sum\limits_{i=1}^{\bar{g}}\boldsymbol{W}_i \quad \notag
\end{equation}
\vspace{-0.05cm}
\begin{displaymath}
\renewcommand{\arraystretch}{2.2}
\qquad\quad=\left\{ \begin{array}{l}
\displaystyle \sum\limits_{i=1}^{\bar{g}}\big (\boldsymbol{W}_{i}^{\delta}+\boldsymbol{W}_{i}^{\omega} \big),  \qquad\qquad\qquad\qquad \textrm{for}\;\textrm{\textbf{M}}_1 \\
\displaystyle \sum\limits_{i=1}^{\bar{g}}\big (\boldsymbol{W}_{i}^{e_R}+\boldsymbol{W}_{i}^{e_I}+\boldsymbol{W}_{i}^{i_R}+\boldsymbol{W}_{i}^{i_I} \big), \;\,\, \textrm{for}\;\textrm{\textbf{M}}_2
\end{array}\right.
\end{displaymath}
where $\boldsymbol{W}^o$ is the empirical observability Gramian for output $o$; $p$ is the number of outputs and $p=2\bar{g}$ for $\textrm{\textbf{M}}_1$  and $p=4\bar{g}$ for $\textrm{\textbf{M}}_2$;
$\boldsymbol{W}_{i}$ is the empirical observability Gramian for all the outputs at generator $i$;
$\boldsymbol{W}_{i}^{\delta}$, $\boldsymbol{W}_{i}^{\omega}$, $\boldsymbol{W}_{i}^{e_R}$, $\boldsymbol{W}_{i}^{e_I}$, $\boldsymbol{W}_{i}^{i_R}$, and $\boldsymbol{W}_{i}^{i_I}$ are respectively
the empirical observability Gramians for the rotor angle, rotor speed, the real or imaginary part of the terminal voltage phaosr, and the real and imaginary part of the terminal current phasor at generator $i$.

Based on this property the optimization problem (\ref{opt}) can be rewritten as
\begin{align}\label{opt1}
&\max\limits_{\boldsymbol{z}} \det \, \sum_{i=1}^g z_i \boldsymbol{W}_{i}(\boldsymbol{z}) \nonumber \\
\;\;\textrm{s.t.}\; & \qquad \sum_{i=1}^g  z_i = \bar{g}  \\
&\qquad z_i \in \{0,1\}, \;\; i=1,\cdots,g \nonumber
\end{align}
where $\boldsymbol{W}_{i}$ is the empirical observability Gramian by only placing one PMU at generator $i$.

The determinant of a matrix is a high-degree polynomial and its absolute value
can be too small or too huge to be represented as a standard double-precision floating-point number.
By contrast, the logarithm of the determinant can be much easier to handle.
Thus we can equivalently rewrite the optimization problem in (\ref{opt1}) as
\begin{align}\label{opt2}
&\min\limits_{\boldsymbol{z}} -\log \det \, \sum_{i=1}^g z_i \boldsymbol{W}_{i}(\boldsymbol{z}) \nonumber \\
\;\;\textrm{s.t.}\; & \qquad \sum_{i=1}^g  z_i = \bar{g}  \\
&\qquad z_i \in \{0,1\}, \;\; i=1,\cdots,g. \nonumber
\end{align}

To summarize, the optimal PMU placement method based on the maximization of the determinant of
the empirical observability Gramian can be implemented in the following two steps.

\begin{enumerate} \renewcommand{\labelitemi}{$\bullet$}
\vspace{0.1cm}
\item \textbf{Calculate empirical observability Gramian} \\
The empirical observability Gramian calculation in (\ref{gd}) is implemented based on emgr (Empirical Gramian Framework) [\ref{emgr}] on time interval $[0,t_f]$. In this paper $t_f$ is chosen as 5 seconds.
$\Delta t_k$ in (\ref{gd}) can take different values according to the required accuracy.
$\boldsymbol{x}_0$ in (\ref{gd}) is chosen as the steady state under typical power flow conditions, which is denoted by $\boldsymbol{x}_0^{ty}$.
We only need to calculate the empirical observability Gramians for placing one PMU at
one of the generators and there is no need to compute all the combinations of PMU placements.
$T^n$ and $M$ that are used to defined the empirical Gramian in (\ref{gd}) are chosen as
\begin{align}
&T^n=\{I_n,-I_n\}  \\
&M=\{0.25,0.5,0.75,1.0\}
\end{align}
where $I_n$ is the identity matrix with dimension $n$ and $n=2\,g$ for model $\textrm{\textbf{M}}_1$ and $n=4\,g_4+2\,g_2$ for model $\textrm{\textbf{M}}_2$. For $T^n$, $I_n$ and $-I_n$ separately correspond to perturbations in the state variables in positive and negative directions. The $M$ chosen here is the default form in emgr, for which the subdivision of the scales of the states is linear and the smallest and biggest perturbation sizes are respectively 0.25 and 1.0.

\vspace{0.1cm}
\item \textbf{Solve MAX-DET optimization problem} \\
MAX-DET problem with continuous variables is convex optimization problem and
can be solved by the interior-point methods [\ref{boyd}] or the Newton-CG primal proximal point algorithm [\ref{wang_det}].
However, the mixed-integer MAX-DET problem in (\ref{opt2}) is nonconvex and cannot be solved by the above-metioned methods.
Therefore, in this paper we resort to the blackbox optimization method and solve (\ref{opt2}) by using the NOMAD solver [\ref{nomad}], which is a derivative-free global mixed integer nonlinear programming solver and is
called by the OPTI toolbox [\ref{opti}].
The NOMAD solver implements the Mesh Adaptive Direct Search (MADS) algorithm [\ref{mads}],
a derivative-free direct search method with a rigorous convergence theory based on the nonsmooth calculus [\ref{clarke}], and aims for the best possible solution with a small number of evaluations.
The advantages of NOMAD and MADS can be summarized as follows:
\begin{enumerate}
\item Under mild hypotheses, the algorithm globally converges to a point satisfying local optimality conditions
based on local properties of the functions defining the problem [\ref{nomad}].
\item Although the MADS algorithm was developed for continuous variables,
the binary control variables in our problem can be easily handled
by using minimal mesh sizes of 1 [\ref{nomad}].
\item NOMAD also includes a Variable Neighborhood Search (VNS) algorithm [\ref{audet}],
which is based on the VNS metaheuristic [\ref{vns}].
This search strategy perturbs the current iterate and conducts poll-like descents
from the perturbed point, allowing an escape from local optima on which the algorithm may be trapped [\ref{nomad}].
\end{enumerate}

Specifically, MADS [\ref{nomad}], [\ref{mads}] is an iterative method where the objective function and constraints are evaluated at some trial points lying on a mesh whose discrete structure is defined at iteration $k$ by
\begin{equation}
M_k=\mycup_{z\in V_k}\{z+\Delta_k^m D u: u\in \mathbb{N}^{n_D}\}
\end{equation}
where $\Delta_k^m \in \mathbb{R}^+$ is the mesh size parameter, $V_k$ is the set of points where the objective function and constraints have been evaluated by the start of iteration $k$, $V_0$ contains the starting points,
and $D=G\,U$ is the set of $n_D$ mesh directions, which is the product of some fixed nonsingular generating matrix $G$ and an integer vector $U$.

Each MADS iteration is composed of three steps: the poll, the search, and updates.
The search step is flexible and allows the creation of trial points anywhere on the mesh.
The poll step is more rigidly defined since the convergence analysis relies on it.
It explores the mesh near the current iterate $z_k$ with the following set of poll trial points:
\begin{equation}
P_k=\{z_k+\Delta_k^m d: d\in D_k\} \subset M_k
\end{equation}
where $D_k$ is the set of poll directions, each column of which is an integer combination of the columns of $D$.
Points of $P_k$ are generated so that their distance to the poll center $x_k$ is bounded below by the poll size
parameter $\Delta_k^p\in \mathbb{R}^+$.

At the end of iteration $k$, an update step determines the iteration status (success
or failure) and the next iterate $x_{k+1}$ is chosen. It corresponds to the most promising
success or stays at $x_k$. The mesh size parameter is also updated with
\begin{equation}
\Delta_{k+1}^m=\tau^{\alpha_k}\Delta_k^m
\end{equation}
where $\tau>1$ is a fixed rational number and $\alpha_k$ is a finite integer, which is positive or null if iteration $k$
succeeds and strictly negative if the iteration fails.
The poll size parameter is also updated in this step according to rules depending on the
implementation.

A high-level description of the algorithm can be given as follows [\ref{nomad}].
  \begin{algorithm}[H]
  \renewcommand\thealgorithm{}
   \caption{MADS} \label{MADS}
    \begin{algorithmic}
      \STATE -- \sc{Initialization}: \textnormal{Let $z_0 \in V_0$ be a starting point,
      $G$ and $U$ be the matrices used to define $D$,
      and $\tau$ be the rational number used to update the mesh size parameter.
      Set $\Delta_0^m$, $\Delta_0^p$, and the iteration counter $k\leftarrow 0$.}
      \vspace{0.1cm}
      \STATE -- Search and poll: \textnormal{Perform the search and poll steps (or only part of them) until an improved mesh point $x_{k+1}$ is found on the mesh $M_k$ or until all trial points are visited.}
      \STATE -- Updates: \textnormal{Update $\Delta_{k+1}^m$, $\Delta_{k+1}^p$, and $V_{k+1}$. \\
      \;Set $k \leftarrow k+1$ and go back to the search and poll steps.}
    \end{algorithmic}
  \end{algorithm}

\end{enumerate}

\subsection{Validation} \label{vali}

The proposed optimal PMU placement method can be validated by
performing dynamic state estimation with SR-UKF [\ref{sr_ukf}].
The following two methods are used to generate a dynamic response.
Method 1 perturbs some of the initial angles and is only used for demonstration while
Method 2 applies a three-phase fault and is thus more realistic.

\begin{enumerate} \renewcommand{\labelitemi}{$\bullet$}
\vspace{0.1cm}
\item \textbf{Method 1--Perturbing angles} \\
The initial states are perturbed by changing some randomly selected angles in the following way
\begin{equation}
\delta_i^\textrm{pert}=\delta_i^{0}+e
\end{equation}
where $i$ is randomly selected among all generators, $\delta_i^{\textrm{pert}}$ is the perturbed angle,
$\delta_i^{0}$ is the steady state angle, and $e$ follows uniform distribution as
\begin{equation}
e\sim U(-|\delta_i^0|,|\delta_i^0|).
\end{equation}

\vspace{0.1cm}
\item \textbf{Method 2--Applying faults} \\
A three-phase fault is applied on a line at one end and is cleared at near and remote end after $0.05s$ and $0.1s$. Here, we do not consider the fault on lines either bus of which is a generator terminal bus because this can lead to the tripping of a generator.
Dynamic state estimation is then performed by SR-UKF on the post-contingency system.
For SR-UKF the mean of the system states at time step 0 is set to be the pre-contingency states,
which can be quite different from the real system states, thus making the dynamic state estimation
very challenging.
\end{enumerate}

To compare the estimation results we define the following system state error
\addtolength{\jot}{1em}
\begin{align}
e_x = \sqrt{\frac{\sum\limits_{i=1}^g\sum\limits_{t=1}^{T_s}(x_{i,t}^{\textrm{est}}-x_{i,t}^{\textrm{true}})^2}{g\,T_s}}
\end{align}
where $x$ is one type of states and can be $\delta$, $\omega$, $e'_q$, or $e'_d$;
$x_{i,t}^{\textrm{est}}$ is the estimated state and $x_{i,t}^{\textrm{true}}$ is the corresponding true value
for the $i$th generator at time step $t$; $T_s$ is the total number of time steps.

We also count the number of convergent states ($\delta$, $\omega$, $e'_q$, or $e'_d$),
which is defined as the states whose differences between the estimated and the true states
in the last 1 second are less than $\epsilon\%$ of the absolute value of the true states.
Obviously the greater the number of convergent states, the better the state estimation result is.

Since where the perturbation is applied and the size of the perturbation are random,
dynamic state estimation is performed for a large number of times in order to get a reliable
conclusion for the comparison of estimation results under different PMU placements.
Note that for different number of PMUs the list of perturbed generators and the perturbation size are the same when using Method 1 and the locations to apply faults and the time to clear faults are the same when using Method 2.

\subsection{Robustness} \label{robustness}

By using the method proposed in Sections \ref{formulate}--\ref{impl} we can obtain the optimal PMU placement for placing $\bar{g}=1,\cdots,g-1$ PMUs under typical power flow conditions. Denote this PMU placement by $\textrm{\textbf{OPP}}^{ty}_{\bar{g}}$, which is the set of generators where PMUs are installed.
However, after a long time the system can significantly change. For example, the loads at some buses might greatly change and even new generators or new transmission lines might be built.
These changes will lead to the change of the unperturbed initial state $\boldsymbol{x}_0$ or even the system dynamics and will further influence the calculation of the empirical observability Gramian and the obtained optimal PMU placement.
In this case the $\textrm{\textbf{OPP}}^{ty}_{\bar{g}}$ obtained based on the current system states
might not be able to make the system state well observed. In order to solve this problem a new optimal PMU placement can be obtained by solving the following optimization problem
\begin{align}\label{opt3}
&\min\limits_{\boldsymbol{z}} -\log \det \, \sum_{i=1}^{g+g_a} z_i \boldsymbol{W}_{o,i}(\boldsymbol{z}) \nonumber \\
\;\;\textrm{s.t.}\; & \qquad \sum_{i=1}^{g+g_a}  z_i = \bar{g}+\bar{g}_a \\
&\qquad z_i=1, \qquad\;\,  i \in \textrm{\textbf{OPP}}^{ty}_{\bar{g}}  \nonumber \\
&\qquad z_i \in \{0,1\}, \;\; i \notin \textrm{\textbf{OPP}}^{ty}_{\bar{g}} \nonumber
\end{align}
where $g_a$ is the number of newly added generators and $\bar{g}_a$ is the number of additional PMUs to be installed.
The optimization problem (\ref{opt3}) keeps the existing PMUs and try to find the best $\bar{g}_a$ locations to install new PMUs in order to maximize the observability of the new system.

However, in this section we mainly discuss the robustness of the proposed optimal PMU placement method when the system does not significantly change.
Specifically we explore the robustness of the proposed optimal PMU placement method under load fluctuations and contingencies, which respectively correspond to small disturbances and big disturbances.
Under these disturbances the unperturbed state $\boldsymbol{x}_0$ in (\ref{gd})
can vary from the steady state under typical power flow conditions $\boldsymbol{x}_0^{ty}$
and when there are contingencies even the system dynamics can change since the topology of the post-contingency system can be different from the pre-contingency system. The robustness of the optimal PMU placement means that the obtained optimal PMU placement will be almost unchanged under these disturbances.

\begin{enumerate} \renewcommand{\labelitemi}{$\bullet$}
\vspace{0.1cm}
\item \textbf{Small disturbance--Load fluctuations} \\
The fluctuations in the loads can be achieved by multiplying the loads under typical power flow conditions by a factor [\ref{OPA}], [\ref{bp13}]
\begin{subequations}
\begin{align}
P_i^{fluc} &= \alpha_i \,P_i^{ty} \\
Q_i^{fluc} &= \alpha_i \,Q_i^{ty}
\end{align}
\end{subequations}
where $i$ belongs to the set of load buses, $P_i^{ty}$ and $Q_i^{ty}$ are real and reactive loads of bus $i$ under typical power flow conditions, $P_i^{fluc}$ and $Q_i^{fluc}$ are real and reactive loads of bus $i$ with fluctuations, and $\alpha_i$ is uniformly distributed in $[2-\gamma,\gamma]$.
We choose $\gamma$ as 1.05 and the real and reactive loads of all load buses will uniformly fluctuate between 95\% and 105\% of their typical power flow loads.

The steady state under load fluctuations $\boldsymbol{x}_0^{fluc}$ is used as
the unperturbed state $\boldsymbol{x}_0$ in (\ref{gd}).

\vspace{0.1cm}
\item \textbf{Big distrubance--Contingencies} \\
Similar to Section \ref{vali}, a three-phase fault is applied on a line at one end and is cleared at near and remote end after $0.05s$ and $0.1s$ and
the faults are also not applied on lines either bus of which is generator terminal buses to avoid the tripping of generators.
The system state after the fault is cleared will be different from the steady state $\boldsymbol{x}_0^{ty}$.
This state is used as the initial unperturbed state $\boldsymbol{x}_0$ in (\ref{gd}) and is denoted by $\boldsymbol{x}_0^{cont}$.
Moreover, the admittance matrix $\boldsymbol{\overline{Y}}$ for the post-contingency system can be different from the pre-contingency system and will further change the system dynamics in both (\ref{gen model1}) and (\ref{gen model}).

\end{enumerate}

The observability Gramian in (\ref{gd}) can be calculated for typical power flow conditions without fluctuations or contingencies and for the above two cases with fluctuations or contingencies and is then used to determine the optimal PMU placement by solving the optimization problem in (\ref{opt2}). The optimal PMU placement for placing $\bar{g}$ PMUs in the two cases with disturbances are separately denoted by $\textrm{\textbf{OPP}}^{fluc}_{\bar{g}}$ and $\textrm{\textbf{OPP}}^{cont}_{\bar{g}}$.
The robustness of the proposed optimal PMU placement method can be verified by comparing $\textrm{\textbf{OPP}}^{fluc}_{\bar{g}}$ and $\textrm{\textbf{OPP}}^{cont}_{\bar{g}}$ with $\textrm{\textbf{OPP}}^{ty}_{\bar{g}}$.

Specifically, the ratios of PMUs placing at the same locations between the typical power flow case and the load fluctuation and contingency case are
\begin{subequations} {\label{ratios}}
\begin{align}
\mathcal{R}_{\bar{g}}^{fluc}&=\frac{\textrm{card}(\textrm{\textbf{OPP}}^{ty}_{\bar{g}} \cap \textrm{\textbf{OPP}}^{fluc}_{\bar{g}})}{\bar{g}} \\
\mathcal{R}_{\bar{g}}^{cont}&=\frac{\textrm{card}(\textrm{\textbf{OPP}}^{ty}_{\bar{g}} \cap \textrm{\textbf{OPP}}^{cont}_{\bar{g}})}{\bar{g}}
\end{align}
\end{subequations}
where $\textrm{card}(\cdot)$ denotes the cardinality of a set, which is a measure of the number of elements of the set. If the ratios defined in (\ref{ratios}) are close to 1.0, $\textrm{\textbf{OPP}}^{ty}_{\bar{g}}$ is close to $\textrm{\textbf{OPP}}^{fluc}_{\bar{g}}$ or $\textrm{\textbf{OPP}}^{cont}_{\bar{g}}$.

Besides, the robustness of the proposed optimal PMU placement method can also be verified
by comparing the logarithm of the determinant of the empirical observability Gramian obtained from the two cases with disturbances for $\textrm{\textbf{OPP}}^{ty}_{\bar{g}}$ and $\textrm{\textbf{OPP}}^{fluc}_{\bar{g}}$ or $\textrm{\textbf{OPP}}^{ty}_{\bar{g}}$ and $\textrm{\textbf{OPP}}^{cont}_{\bar{g}}$. Specifically, the corresponding $\boldsymbol{z}$ in (\ref{opt2}) (separately denoted by $\boldsymbol{z}^{ty}_{\bar{g}}$, $\boldsymbol{z}^{fluc}_{\bar{g}}$, and $\boldsymbol{z}^{cont}_{\bar{g}}$) can be obtained from $\textrm{\textbf{OPP}}^{ty}_{\bar{g}}$, $\textrm{\textbf{OPP}}^{fluc}_{\bar{g}}$, or $\textrm{\textbf{OPP}}^{cont}_{\bar{g}}$ since $z_i=1$ if $i \in \textrm{\textbf{OPP}}_{\bar{g}}$ and $z_i=0$ otherwise, where $\textrm{\textbf{OPP}}_{\bar{g}}$ can be the optimal placement in any of the three cases.
Then the objective function in (\ref{opt2}), for which $\boldsymbol{W}_{o,i}(\boldsymbol{z})$ is the empirical observability Gramian for the small or big disturbance case, can be evaluated for $\boldsymbol{z}^{ty}_{\bar{g}}$ and $\boldsymbol{z}^{fluc}_{\bar{g}}$ or $\boldsymbol{z}^{ty}_{\bar{g}}$ and $\boldsymbol{z}^{cont}_{\bar{g}}$. The logarithm of the determinant of the empirical observability Gramian can finally be obtained  as the opposite of the optimal value in (\ref{opt2}). Here robustness means that the logarithms of the determinant of the empirical observability Gramian are similar for $\textrm{\textbf{OPP}}^{ty}_{\bar{g}}$ and $\textrm{\textbf{OPP}}^{fluc}_{\bar{g}}$ in the small disturbance case or for $\textrm{\textbf{OPP}}^{ty}_{\bar{g}}$ and $\textrm{\textbf{OPP}}^{cont}_{\bar{g}}$ in the big disturbance case.

\subsection{Comparison with Existing Methods}

As is mentioned in Section \ref{intro}, there are few methods on PMU placement for dynamic state estimation.
In this section we will briefly compare the proposed method in this paper with some existing methods.

In [\ref{kamwa}] numerical PMU configuration algorithms are proposed to maximize the overall sensor response
while minimizing the correlation among sensor outputs so as to minimize the redundant information provided by multiple sensors. The optimal PMU placement problem is not tackled directly but is solved by an sequential ``greedy" heuristic algorithm [\ref{greedy}].

In [\ref{zhang1}] and [\ref{zhang2}] the performance of multiple optimal PMU placements obtained by mixed integer programming [\ref{ob2}], [\ref{ob3}] are evaluated by using a stochastic estimate of the steady-state error covariance. They only compare some already obtained PMU placements but do not propose a method to determine the optimal PMU placement.

In [\ref{huang}] a PMU placement method is proposed to ensure a satisfactory state tracking performance.
It depends on a specific Kalman filter and tries to find a PMU placement strategy with small tracking error. However, it does not answer the question why some specific PMU placement can guarantee small tracking error. The optimal PMU placement problem is solved by a sequential heuristic algorithm.

By contrast, the advantages of the optimal PMU placement method proposed in this paper can be summarized as follows:
\begin{enumerate}
\item The proposed method has a quantitative measure of observability, the determinant of the empirical observability Gramian, which makes it possible to optimize PMU locations from the point view of the observability of nonlinear systems.
\item The proposed method efficiently solves the PMU placement problem by the NOMAD solver [\ref{nomad}], which is a derivative-free global mixed integer nonlinear programming solver. This is significantly different from similar work solely based on the sequential heuristic algorithm [\ref{kamwa}], [\ref{huang}].
\item The proposed method does not need to consider different load levels or contingencies but only needs to deal with the system under typical power flow conditions. This is true especially after the robustness of the method discussed in Section \ref{robustness} is well validated, for which the results will be given in Section \ref{robu}. This is different from [\ref{kamwa}] for which the PMU placement is obtained by using the system response under many contingencies.
\item The proposed method does not need to perform dynamic state estimation and thus does not depend on the specific realization of Kalman filter. Thus the obtained optimal PMU placement can be applied to any type of Kalman filter, which is an important advantage over the method in [\ref{huang}].
\end{enumerate}

Even compared with more general sensor placement methods with applications other than power system dynamic state estimation, the proposed method in this paper also has obvious advantages.
\begin{enumerate}
\item In [\ref{singh1}] the trace of the empirical observability Gramian is used as measure of observability and the objective function of the optimal sensor placement problem, which is used to find optimal sensor placement for the distillation column model. However, as is discussed in Section \ref{formulate}, the trace cannot tell the existence of a zero eigenvalue. The proposed method can avoid this by using the determinant as the measure of observability.
Moreover, in [\ref{singh1}] the optimization problem is solved by the genetic algorithm, which is often considered as being not time-efficient.
\item In [\ref{det}] the optimal sensor placement for distillation column and packed bed reactor is obtained by maximizing the determinant of the empirical observability Gramian. However, the BONMIN slover [\ref{algo}]
that is chosen to solve the MAX-DET problem can only solve convex mixed integer nonlinear programming. It is not suitable for solving the nonconvex mixed integer MAX-DET problem and can only find local solutions to nonconvex problems. By contrast, the NOMAD solver [\ref{nomad}] is a derivative-free global mixed integer nonlinear programming solver. The examples that show the advantage of the NOMAD solver over the BONMIN solver on solving nonconvex mixed integer nonlinear programming problems can be found on the website of the OPTI toolbox [\ref{opti}] at \url{http://www.i2c2.aut.ac.nz/Wiki/OPTI/index.php/Probs/MINLP}.
\end{enumerate}

\section{Case Studies} \label{case}

In this section the proposed optimal PMU placement method is tested on WSCC 3-machine 9-bus system [\ref{anderson}] and NPCC 48-machine 140-bus system [\ref{pst}]. The empirical observability Gramian calculation and the SR-UKF are implemented with Matlab and all tests are carried out on a 3.4 GHz Intel(R) Core(TM) i7-3770 based desktop.
When calculating the empirical observability Gramian in (\ref{gd}) $\Delta t_k$ is chosen as $1/30s$ for WSCC system and $1/120s$ for NPCC system since it is more difficult to reliably and stably solve the ordinary differential equations for the much bigger NPCC system.
Dynamic state estimation is performed on $[0,5s]$ by using SR-UKF.
By denoting the error of $\delta$, $\omega$, $e'_q$, $e'_d$, $e_R$, $e_I$, $i_R$, and $i_I$ respectively as $r_\delta$, $r_\omega$, $r_{e'_q}$, $r_{e'_d}$, $r_{e_R}$, $r_{e_I}$, $r_{i_R}$, and $r_{i_I}$,
the SR-UKF is set as follows.
\begin{enumerate}
\item For $\textrm{\textbf{M}}_1$ the initial state covariance is set as
\begin{displaymath}
\renewcommand{\arraystretch}{1.4}
\boldsymbol{P}_{0,1} =
\left[ \begin{array}{cc}
r_\delta^2 \it I_g & \boldsymbol{0}_g \\
\boldsymbol{0}_g & r_\omega^2 I_g
\end{array} \right]
\end{displaymath}
and for $\textrm{\textbf{M}}_2$ the initial state covariance is
\begin{displaymath}
\renewcommand{\arraystretch}{1.4}
\boldsymbol{P}_{0,2} =
\left[ \begin{array}{cccc}
r_\delta^2 I_g & \boldsymbol{0}_g & \boldsymbol{0}_{g,g_4} & \boldsymbol{0}_{g,g_4} \\
\boldsymbol{0}_g & r_\omega^2 I_g & \boldsymbol{0}_{g,g_4} & \boldsymbol{0}_{g,g_4} \\
\boldsymbol{0}_{g_4,g} & \boldsymbol{0}_{g_4,g} & r_{e'_q}^2 I_{g_4} & \boldsymbol{0}_{g_4} \\
\boldsymbol{0}_{g_4,g} & \boldsymbol{0}_{g_4,g} & \boldsymbol{0}_{g_4} & r_{e'_d}^2  I_{g_4}
\end{array} \right]
\end{displaymath}
where $I_g$ is an identity matrix with dimension $g$, $\boldsymbol{0}_g$ is a zero matrix with dimension $g$, and $\boldsymbol{0}_{u,v}$ is a zero matrix with dimension $u\times v$.

\item The covariance for the process noise is $\boldsymbol{Q}=10^{-7} I_{n}$ where $n$ is the number of states.

\item For $\textrm{\textbf{M}}_1$ the covariance for the measurement noise is
\begin{displaymath}
\renewcommand{\arraystretch}{1.4}
\boldsymbol{R}_1 =
\left[ \begin{array}{cc}
r_{\delta}^2 I_{\bar{g}} & \boldsymbol{0}_{\bar{g}} \\
\boldsymbol{0}_{\bar{g}} & r_{\omega}^2 I_{\bar{g}}
\end{array} \right]
\end{displaymath}
and for $\textrm{\textbf{M}}_2$ the covariance for the measurement noise is
\begin{displaymath}
\renewcommand{\arraystretch}{1.4}
\boldsymbol{R}_2 =
\left[ \begin{array}{cccc}
r_{e_R}^2 I_{\bar{g}} & \boldsymbol{0}_{\bar{g}} & \boldsymbol{0}_{\bar{g}} & \boldsymbol{0}_{\bar{g}}\\
\boldsymbol{0}_{\bar{g}} & r_{e_I}^2 I_{\bar{g}} & \boldsymbol{0}_{\bar{g}} & \boldsymbol{0}_{\bar{g}} \\
\boldsymbol{0}_{\bar{g}} & \boldsymbol{0}_{\bar{g}} & r_{i_R}^2 I_{\bar{g}} & \boldsymbol{0}_{\bar{g}} \\
\boldsymbol{0}_{\bar{g}} & \boldsymbol{0}_{\bar{g}} & \boldsymbol{0}_{\bar{g}} & r_{i_I}^2 I_{\bar{g}}
\end{array} \right]
\end{displaymath}
where $\bar{g}$ is the number of PMUs.

\item The PMU sampling rate is set to be 30 frames per second for WSCC 3-machine system and
60 frames per second for NPCC 48-machine system. Note that if the practical sampling rate of PMUs in real systems is smaller than 60 frames per second the effective sampling rate can be increased by the interpolation method [\ref{huang1}] which adds pseudo measurement points between two consecutive measurement samples.

\end{enumerate}

\subsection{WSCC 3-Machine System}

The proposed optimal PMU placement method is applied to WSCC 3-machine 9-bus system to decide where PMUs should be
installed to make the system most observable.
This small system is only used for demonstration and thus the simplified generator and measurement model in Section \ref{sim model} is used.
Method 1 in Section \ref{vali} is applied to generate dynamic response
and specifically one of the angles is perturbed.
We choose Method 1 because WSCC 3-machine system is very small and there are not many different cases
if we apply faults by using Method 2.
$\epsilon$ is set as 2 and $r_\delta$ and $r_\omega$ are chosen as $0.5\,\pi /180$ and $10^{-3}\omega_0$. The optimal PMU placements for placing 1 and 2 PMUs are listed in Table \ref{pmu_3}, in which the time for solving the optimization problem (\ref{opt2}) is also listed.

\begin{figure}[!t]
\centering
\includegraphics[width=3.2in]{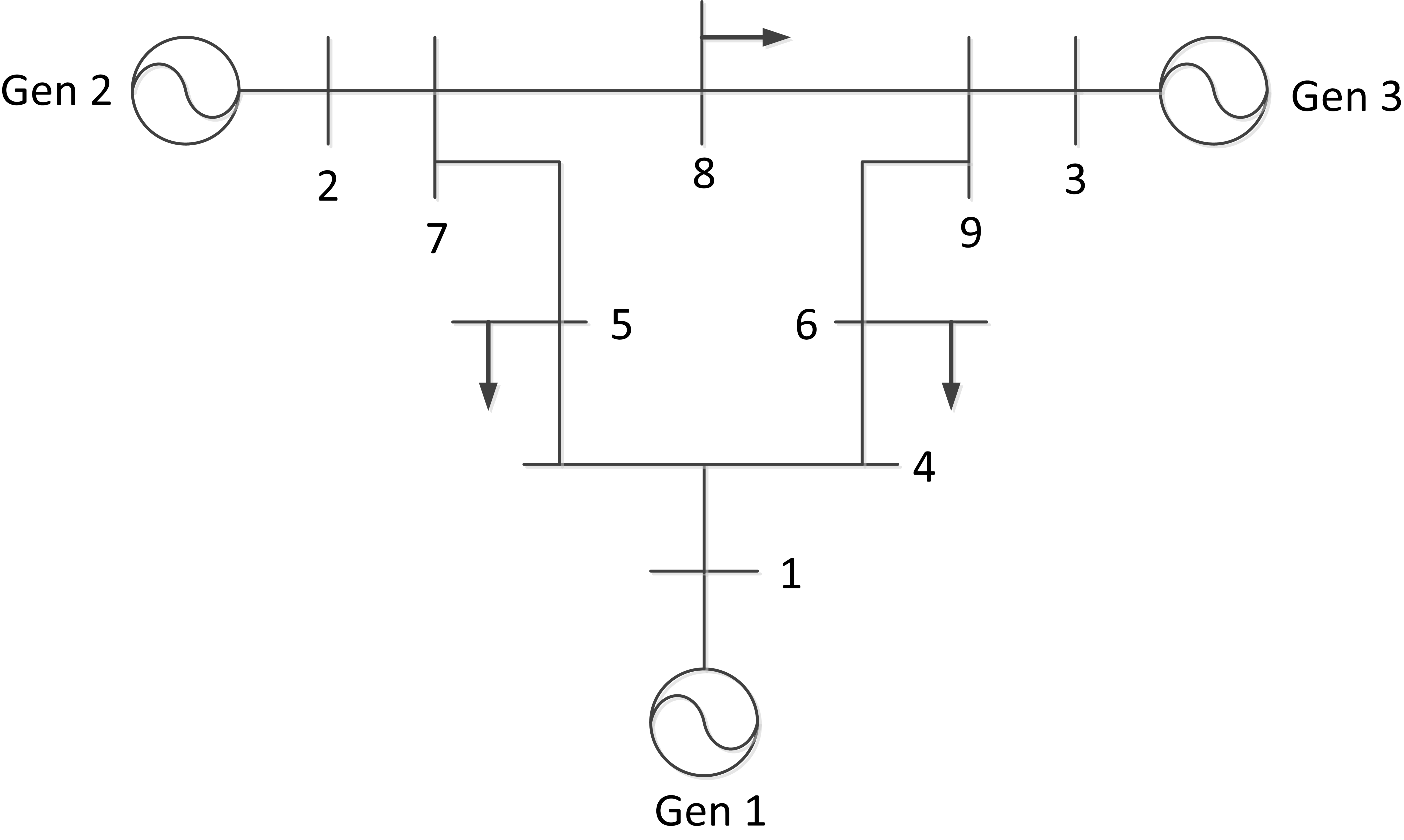}
\caption{WSCC 3-machine 9-bus system.}
\label{wscc3}
\end{figure}

\begin{table}[H]
\renewcommand{\arraystretch}{1.3}
\caption{Optimal PMU Placement for WSCC 3-Machine System}
\label{pmu_3}
\centering
\begin{tabular}{ccc}
\hline
number of PMUs & \tabincell{c}{optimal placement} & time ($s$) \\
\hline
1 & 3 & 2.46 \\
2 & 2,\,3 & 0.076 \\
\hline
\end{tabular}
\end{table}

The logarithms of determinant of the empirical observability Gramian
under different PMU placements are listed in Table \ref{logdet}, in which
the logarithm of determinants under the optimal PMU placement
achieve maximum among those with the same number of PMUs.
We also list the largest and smallest eigenvalue $\sigma_{\textrm{max}}$ and $\sigma_{\textrm{min}}$.
The PMU placement with the greatest logarithm of determinant also corresponds to the greatest $\sigma_{\textrm{max}}$ and $\sigma_{\textrm{min}}$.

\begin{table}[H]
\renewcommand{\arraystretch}{1.3}
\caption{Logarithm of Determinant of the Empirical Observability Gramian under Different PMU Placements for WSCC 3-Machine System}
\label{logdet}
\centering
\begin{tabular}{cccc}
\hline
PMU placement & $\log \det$ & $\sigma_{\textrm{max}}$ & $\sigma_{\textrm{min}}$ \\
\hline
1 & 8.54  & $1.14\times 10^3$ & 0.0082 \\
2 & 19.61 & $1.16\times 10^3$ & 0.43 \\
3 & 22.33 & $1.23\times 10^3$ & 0.57 \\
1,\,2 & 21.34 & $2.30\times 10^3$ & 0.44 \\
1,\,3 & 24.40 & $2.37\times 10^3$ & 0.82 \\
2,\,3 & 26.47 & $2.40\times 10^3$ & 2.15 \\
\hline
\end{tabular}
\end{table}

From Table \ref{logdet} it is seen that the logarithms of determinant
of the empirical observability Gramian for placing one PMU at generator 2 and generator 3
are much greater than that for placing PMU at generator 1, indicating that
placing one PMU at generator 2 or generator 3 should make the system more observable than
placing one PMU at generator 1.

\begin{figure}[!t]
\centering
\includegraphics[width=3.5in]{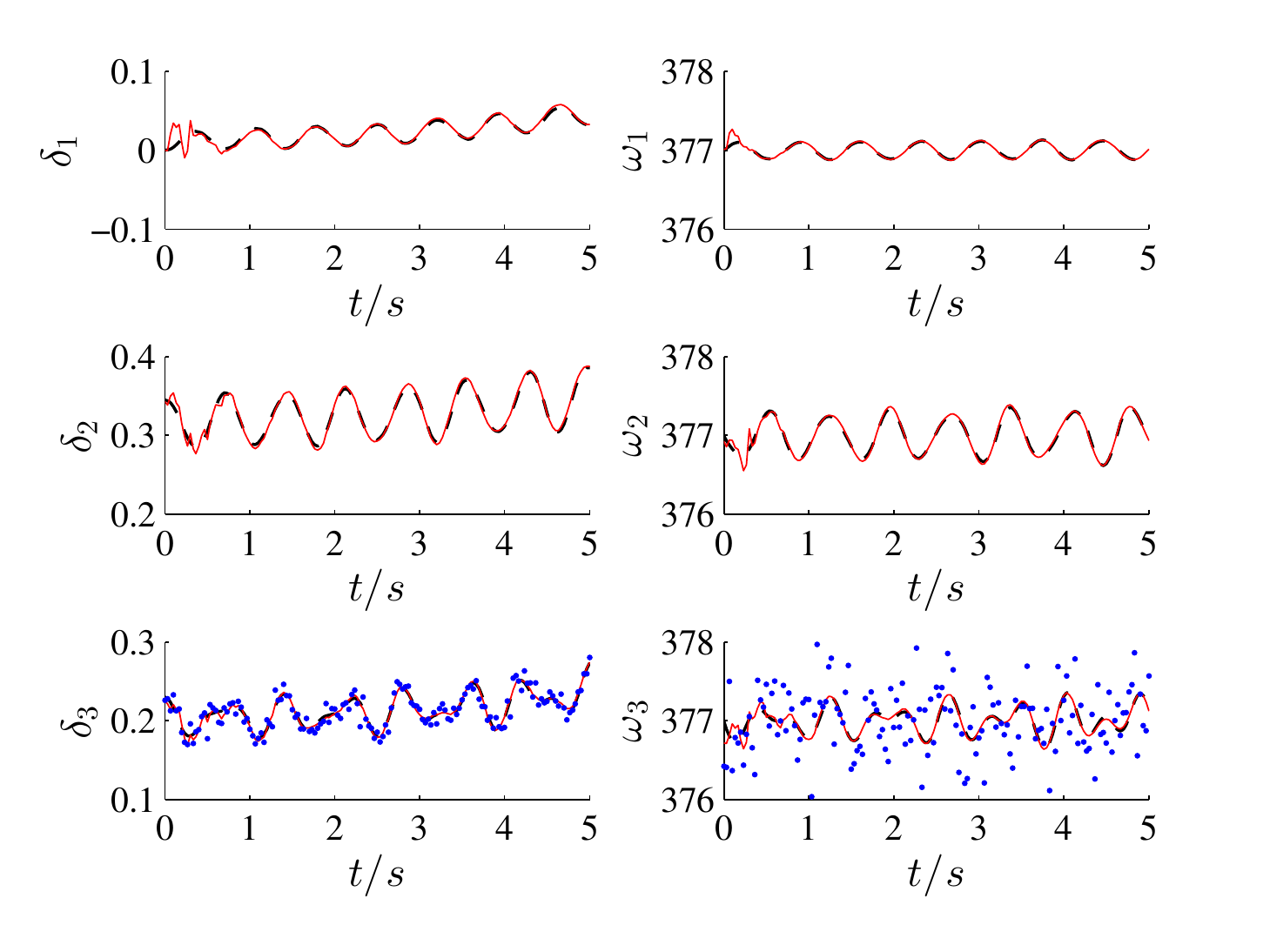}
\caption{Estimation results for placing one PMU at generator 3 for WSCC 3-machine system.}
\label{pmu3}
\end{figure}

\begin{figure}[!t]
\centering
\includegraphics[width=3.5in]{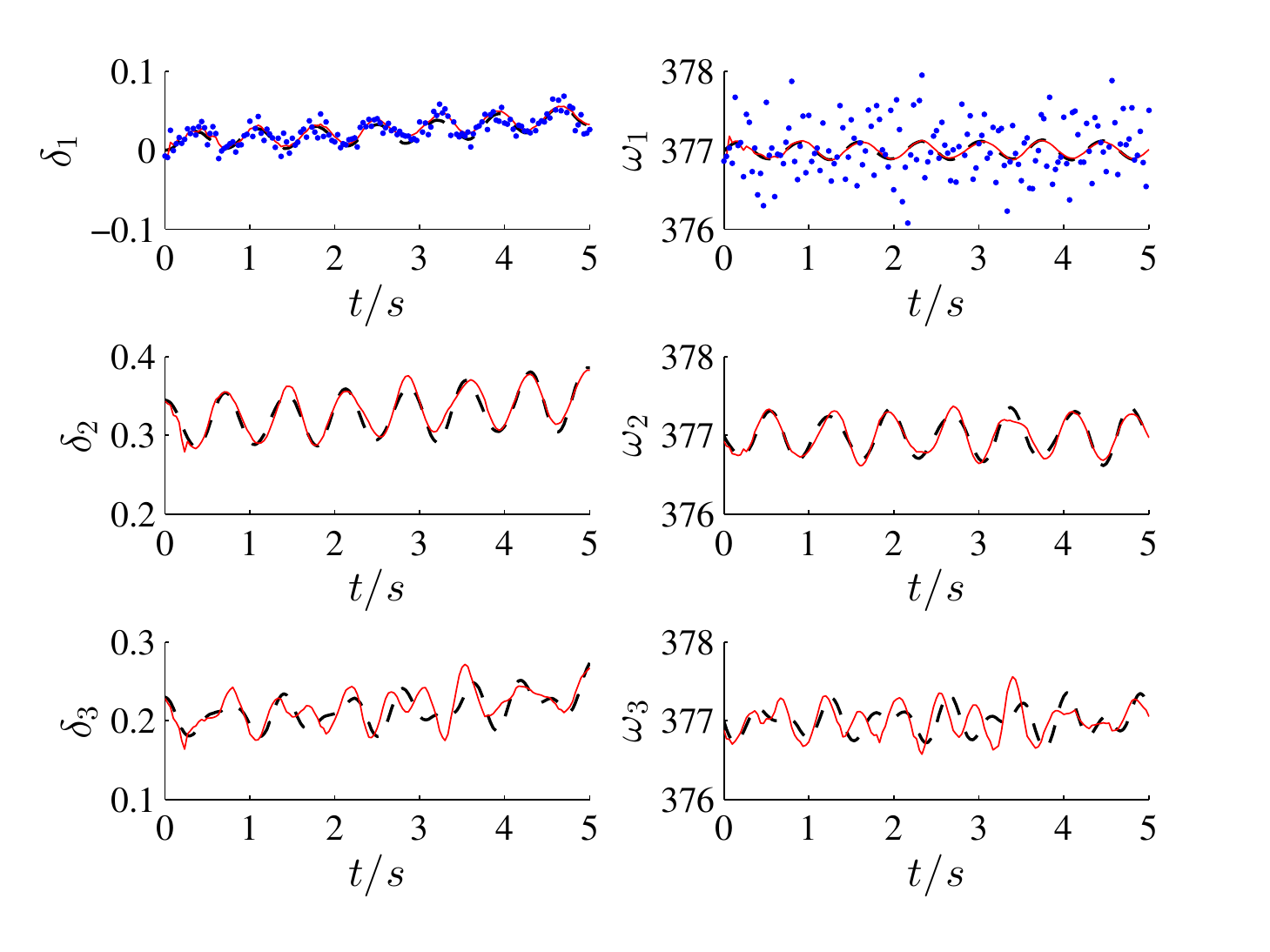}
\caption{Estimation results for placing one PMU at generator 1 for WSCC 3-machine system.}
\label{pmu1}
\end{figure}

To show the difference for placing PMU at generator 3 and generator 1 we present
estimation results for perturbing the angle of generator 1 by decreasing it by 100\%
in Figs. \ref{pmu3} and \ref{pmu1}, in which black dash lines denote real states, red solid lines
denote estimated states, and blue dots denotes the measurements from PMUs.
When placing PMU at generator 3 all the rotor angles and rotor speeds can
converge to the true states quickly while for placing PMU at generator 1 the rotor angles and rotor speeds
where PMUs are not installed are difficult to converge. The rotor angle and rotor speed errors for placing PMU at generator 3 are 0.0044 and 0.041 while for placing PMU at generator 1 are 0.011 and 0.12.

The average value of the rotor angle and rotor speed error $\bar{e}_\delta$ and $\bar{e}_\omega$ and the average value of the number of convergent angles $\overline{N}_\delta$ for performing dynamic state estimation
using SR-UKF for 50 times under the optimal and other PMU placements are listed in Table \ref{error3}.
For each dynamic state estimation one randomly selected angle is perturbed by Method 1 in Section \ref{vali}
in order to generate different cases.
Under optimal PMU placements the rotor angle and rotor speed errors are the smallest and the numbers of convergent angles are the greatest.

\begin{table}[!t]
\renewcommand{\arraystretch}{1.3}
\caption{Error and Number of Convergent Angles under Different PMU Placements for WSCC 3-Machine System}
\label{error3}
\centering
\begin{tabular}{ccccc}
\hline
PMU placement & $\bar{e}_\delta$ & $\bar{e}_\omega$ & $\overline{N}_\delta$ \\
\hline
 1 & 0.015 & 0.19 & 0.88 \\
 2 & 0.0078 & 0.083 & 1.96 \\
 3 & 0.0058 & 0.055 & 2.04 \\
 1,\,2 & 0.0056 & 0.065 & 2.14 \\
 1,\,3 & 0.0044 & 0.042 & 2.22 \\
 2,\,3 & 0.0037 & 0.036 & 2.28 \\
\hline
\end{tabular}
\end{table}

\subsection{NPCC 48-Machine System}

\begin{figure}[!t]
\centering
\includegraphics[width=3.3in]{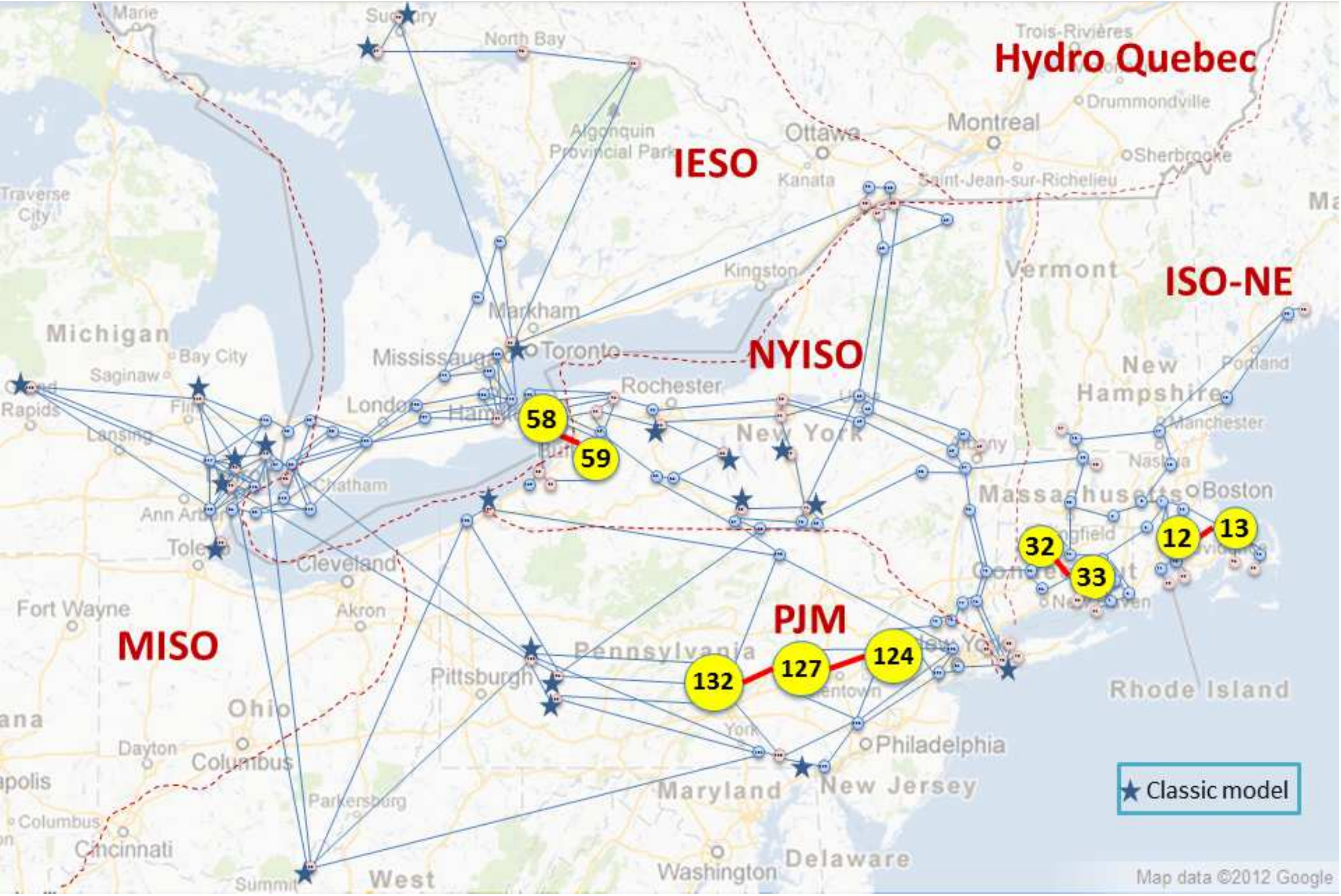}
\caption{Map of the NPCC 48-machine 140-bus system. The stars indicates generators with classical model and the transmission lines with highest line flow and their buses are highlighted.}
\label{npcc48}
\end{figure}

In this section the proposed optimal PMU placement method is tested on NPCC 48-machine system [\ref{pst}], which has 140 buses and
represents the northeast region of the EI system. The map of this system is shown in Fig. \ref{npcc48}.
The realistic generator and measurement model in Section \ref{sim model} is used and
27 generators have fourth-order transient model and the other 21 generators have second-order classical model.
Method 2 in Section \ref{vali} is applied to generate dynamic response.
$\epsilon$ is set to be 2. $r_\delta$ and $r_\omega$ are chosen as $0.5\,\pi /180$ and $10^{-3}\omega_0$. 
$r_{e'_q}$ and $r_{e'_d}$ are set to be $10^{-3}$.
$r_{e_R}$ and $r_{e_I}$ are chosen as $5\times 10^{-2}$ and $r_{i_R}$ and $r_{i_I}$ are chosen as $5\times 10^{-1}$.

The obtained optimal PMU placements are listed in Table \ref{pmu_48}.
For brevity we only list the PMU placements for the number of PMUs between 12 and 24.
We can see that if $\bar{g}_1>\bar{g}_2$ it does not necessarily hold that the optimal PMU placement $\textrm{\textbf{OPP}}^{ty}_{\bar{g}_1} \supset \textrm{\textbf{OPP}}^{ty}_{\bar{g}_2}$.
This is very different from [\ref{kamwa}] and [\ref{huang}], in which a sequential heuristic algorithm is used and the optimal solution can not be guaranteed.

\begin{table}[!t]
\renewcommand{\arraystretch}{1.5}
\caption{Optimal PMU Placement for NPCC 48-Machine System \qquad\qquad (Number of PMUs = $12\sim 24$)}
\label{pmu_48}
\centering
\begin{tabular}{cl}
\hline
number of PMUs & \quad\quad\quad \tabincell{c}{optimal placement} \\
\hline
12 & \tabincell{l}{2,\,3,\,6,\,11,\,13,\,16,\,18,\,21,\,27,\,32,\,33,\,44}  \\
13 & \tabincell{l}{2,\,3,\,6,\,11,\,13,\,17,\,18,\,21,\,27,\,32,\,33,\,37,\,44}  \\
14 & \tabincell{l}{2,\,3,\,6,\,11,\,13,\,17,\,18,\,19,\,22,\,27,\,32,\,33,\,37,\,44}  \\
15 & \tabincell{l}{2,\,3,\,6,\,11,\,13,\,16,\,18,\,19,\,21,\,27,\,28,\,32,\,33,\,38,\,44}  \\
16 & \tabincell{l}{2,\,3,\,6,\,12,\,13,\,16,\,18,\,19,\,22,\,27,\,28,\,32,\,33,\,37,\,44,\,45}  \\
17 & \tabincell{l}{2,\,3,\,6,\,9,\,12,\,13,\,17,\,18,\,19,\,21,\,27,\,28,\,32,\,33,\,37,\,44,\,45}  \\
18 & \tabincell{l}{2,\,3,\,6,\,9,\,11,\,13,\,16,\,18,\,19,\,21,\,27,\,28,\,31,\,32,\,33,\\37,\,44,\,45} \\
19 & \tabincell{l}{2,\,3,\,6,\,9,\,11,\,13,\,14,\,17,\,18,\,19,\,21,\,27,\,28,\,31,\,32,\\ 33,\,38,\,44,\,45} \\
20 & \tabincell{l}{2,\,3,\,6,\,9,\,11,\,13,\,14,\,17,\,18,\,19,\,20,\,21,\,27,\,28,\,31,\\ 32,\,33,\,38,\,44,\,45} \\
21 & \tabincell{l}{1,\,2,\,3,\,6,\,9,\,12,\,13,\,14,\,17,\,18,\,19,\,20,\,21,\,27,\,28,\\ 31,\,32 \,33,\,37,\,44,\,45} \\
22 & \tabincell{l}{1,\,2,\,3,\,6,\,9,\,10,\,11,\,13,\,14,\,17,\,18,\,19,\,20,\,21,\,27,\\ 29,\,31,\,32,\,33,\,37,\,44,\,45} \\
23 & \tabincell{l}{1,\,2,\,3,\,4,\,6,\,9,\,10,\,12\,13,\,14,\,16,\,18,\,19,\,20,\,21,\\ 27,\,28,\,31,\,32,\,33,\,37,\,44,\,48} \\
24 & \tabincell{l}{1,\,2,\,3,\,4,\,6,\,9,\,10,\,12,\,13,\,14,\,16,\,18,\,19,\,20,\,21,\\ 27,\,28,\,31,\,32,\,35,\,36,\,38,\,44,\,45} \\
\hline
\end{tabular}
\end{table}

\begin{figure}[!t]
\centering
\includegraphics[width=2.9in]{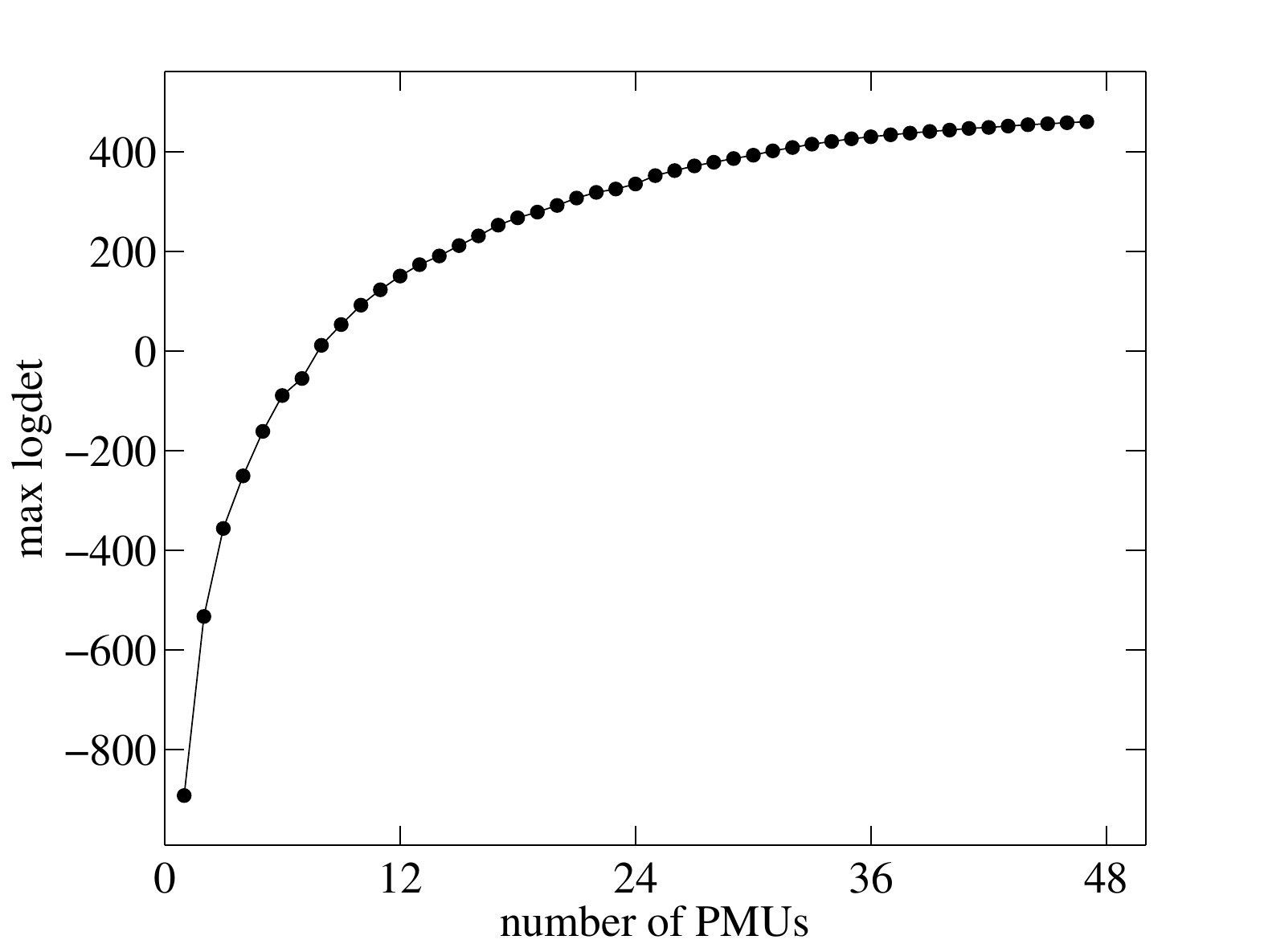}
\caption{Maximum logarithm of determinant of the empirical observability Gramian for different numbers of PMUs for NPCC 48-machine system.}
\label{fval3}
\end{figure}

\begin{figure}[!t]
\centering
\includegraphics[width=2.9in]{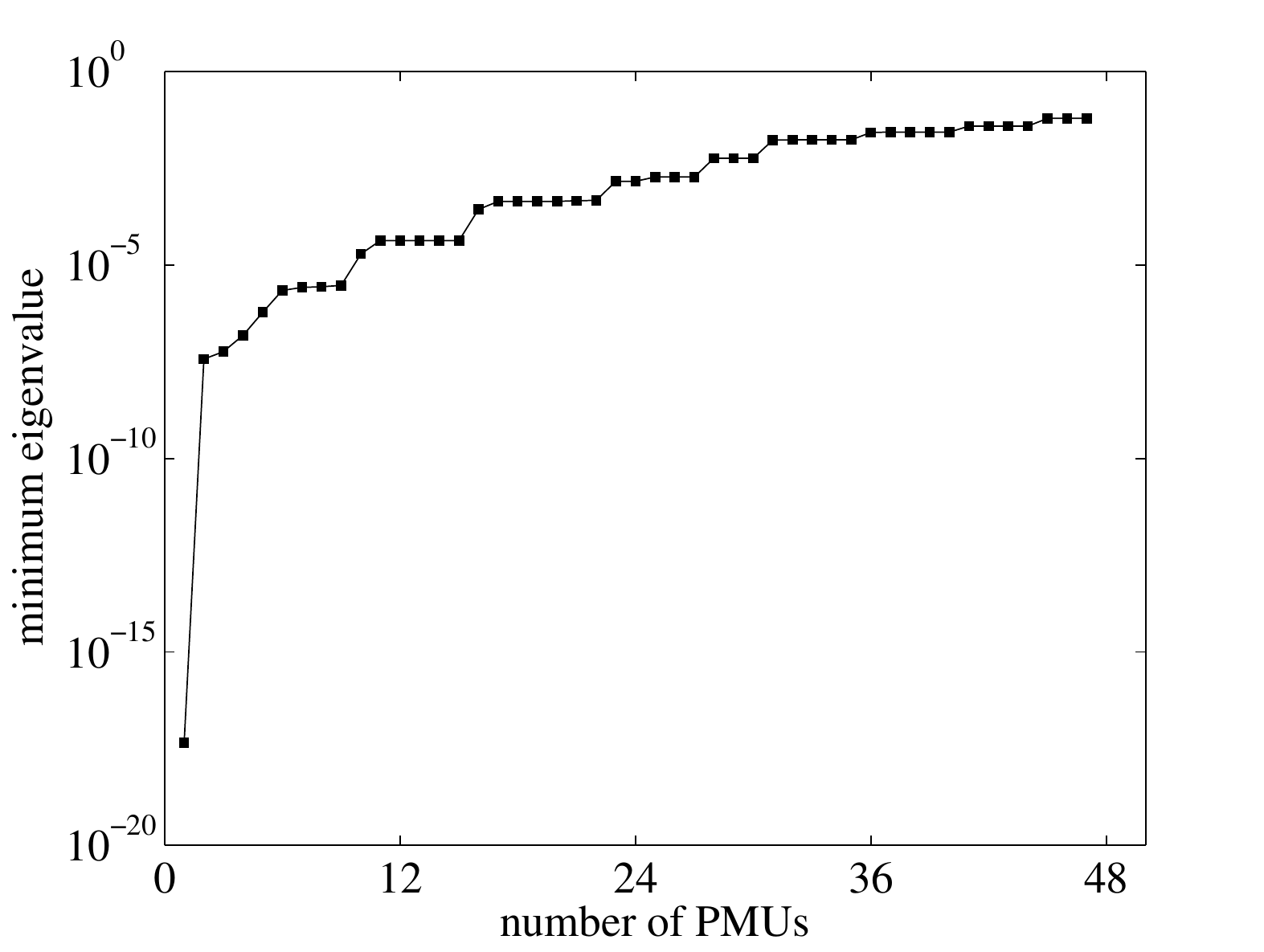}
\caption{Minimum eigenvalue of the empirical observability Gramian for different numbers of PMUs for NPCC 48-machine system.}
\label{minE}
\end{figure}

The maximum logarithms of the determinant of the empirical observability Gramian
for placing different numbers of PMUs $\bar{g}$ are shown in Fig. \ref{fval3}.
The maximum logarithms of the determinant first increase quickly when the number of PMUs is small
and then increases much slower after $\bar{g}$ exceeds around 10.
As is mentioned in Section \ref{formulate}, it is advised that the smallest eigenvalue be verified to be at an acceptable level while using the determinant of Gramian as the measure of observability in order to avoid the situation in which a Gramian has an almost zero eigvenvalue that makes the system practically unobservable.
Therefore, the minimum eigenvalues of the empirical observability Gramian for different number of PMUs are shown in Fig. \ref{minE}, from which we can see that the minimum eigenvalues of the observability Gramian gradually increase with the increase of $\bar{g}$. Except for the cases in which $\bar{g}$ is very small the minimum eigenvalues for most cases can stay at an acceptable level. Specifically, when $\bar{g}=10$, the minimum eigenvalue becomes greater than $10^{-5}$; when $\bar{g}=23$, the minimum eigenvalue becomes greater than $10^{-3}$; finally when $\bar{g} \ge 45$ the minimum eigenvalues are around 0.062.

\begin{figure}[!t]
\centering
\includegraphics[width=2.9in]{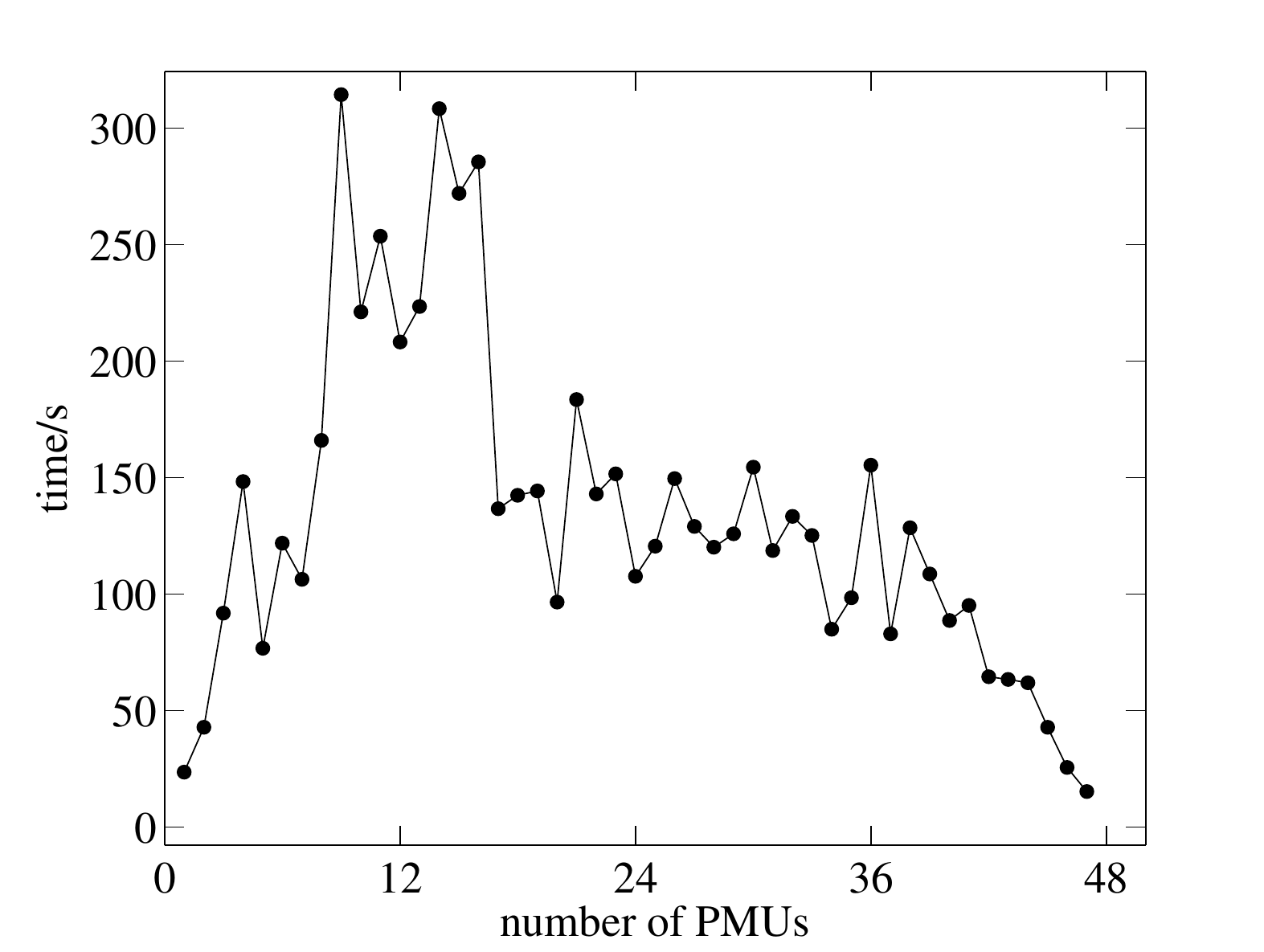}
\caption{Time for solving the PMU placement problem (\ref{opt2}) for different numbers of PMUs for NPCC 48-machine system.}
\label{time}
\end{figure}

The time for calculating the optimization problem (\ref{opt2})
for different numbers of PMUs is shown in Fig. \ref{time}.
Note that for most cases in which $\bar{g}=1 \sim 47$ the number of possible placements are huge and
an exhaustive enumeration of the solution space is unfeasible.
By contrast, by using the NOMAD solver [\ref{nomad}] we are able to solve (\ref{opt2}) very efficiently.
The longest calculation time (314 seconds) corresponds to
placing 9 PMUs and for placing $\bar{g} = 1 \sim 47$ PMUs 89\% of the calculation can be finished in 4 minutes, indicating that the calculation is very time-efficient.

We compare the estimation error of the states and the number of convergent states between the obtained optimal PMU placements and random PMU placements for which the same number of PMUs are placed at randomly selected generators. Method 2 in Section \ref{vali} is applied to generate dynamic response.
50 cases are created and for each of them a three-phase fault is applied at the from bus of one of the 50
branches with highest line flows and is cleared at near and remote end after $0.05s$ and $0.1s$.

\begin{figure}[!t]
\centering
\includegraphics[width=2.9in]{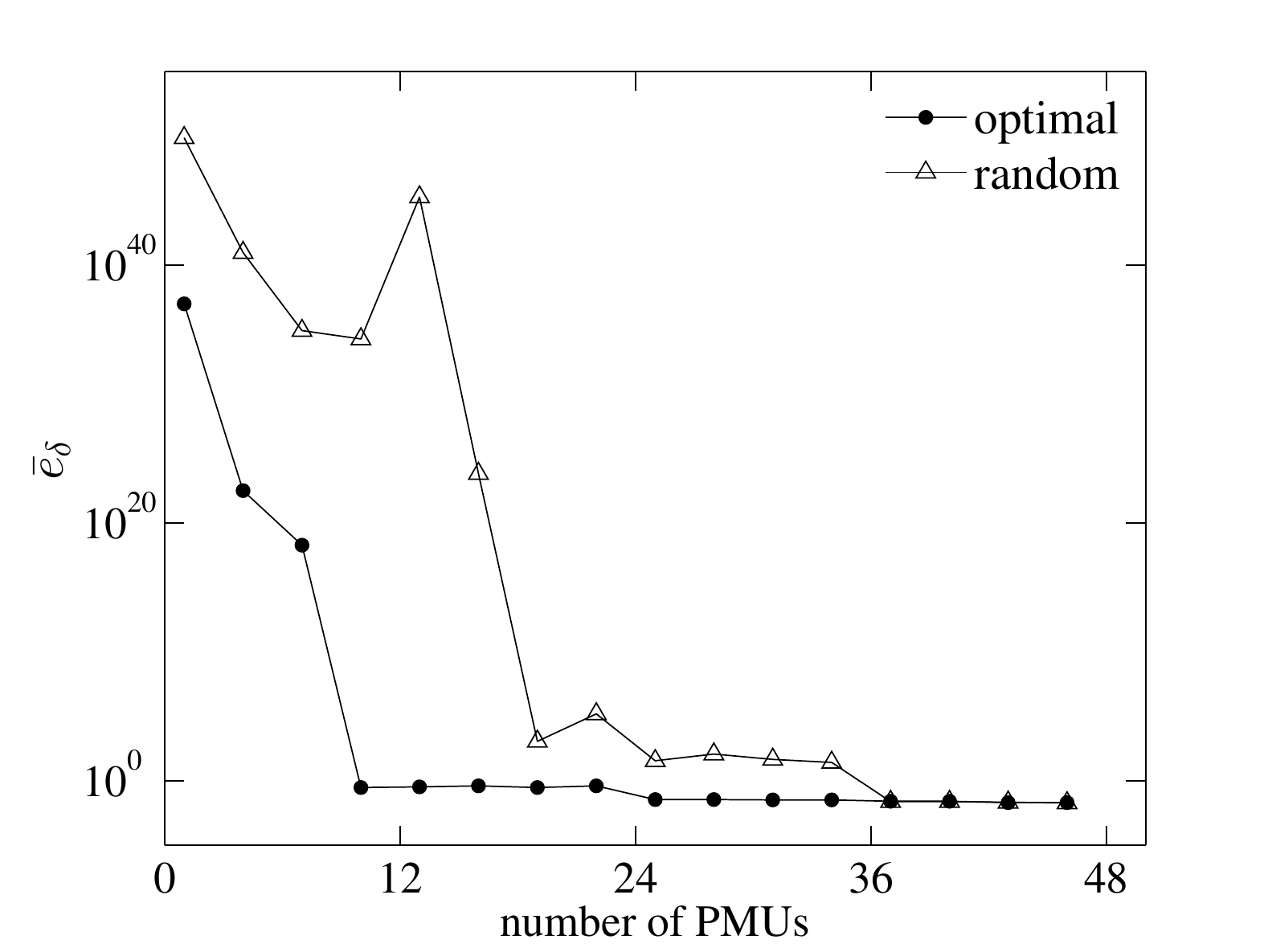}
\caption{Error of rotor angles under optimal and random PMU placements for NPCC 48-machine system.}
\label{d}
\end{figure}

\begin{figure}[!t]
\centering
\includegraphics[width=2.9in]{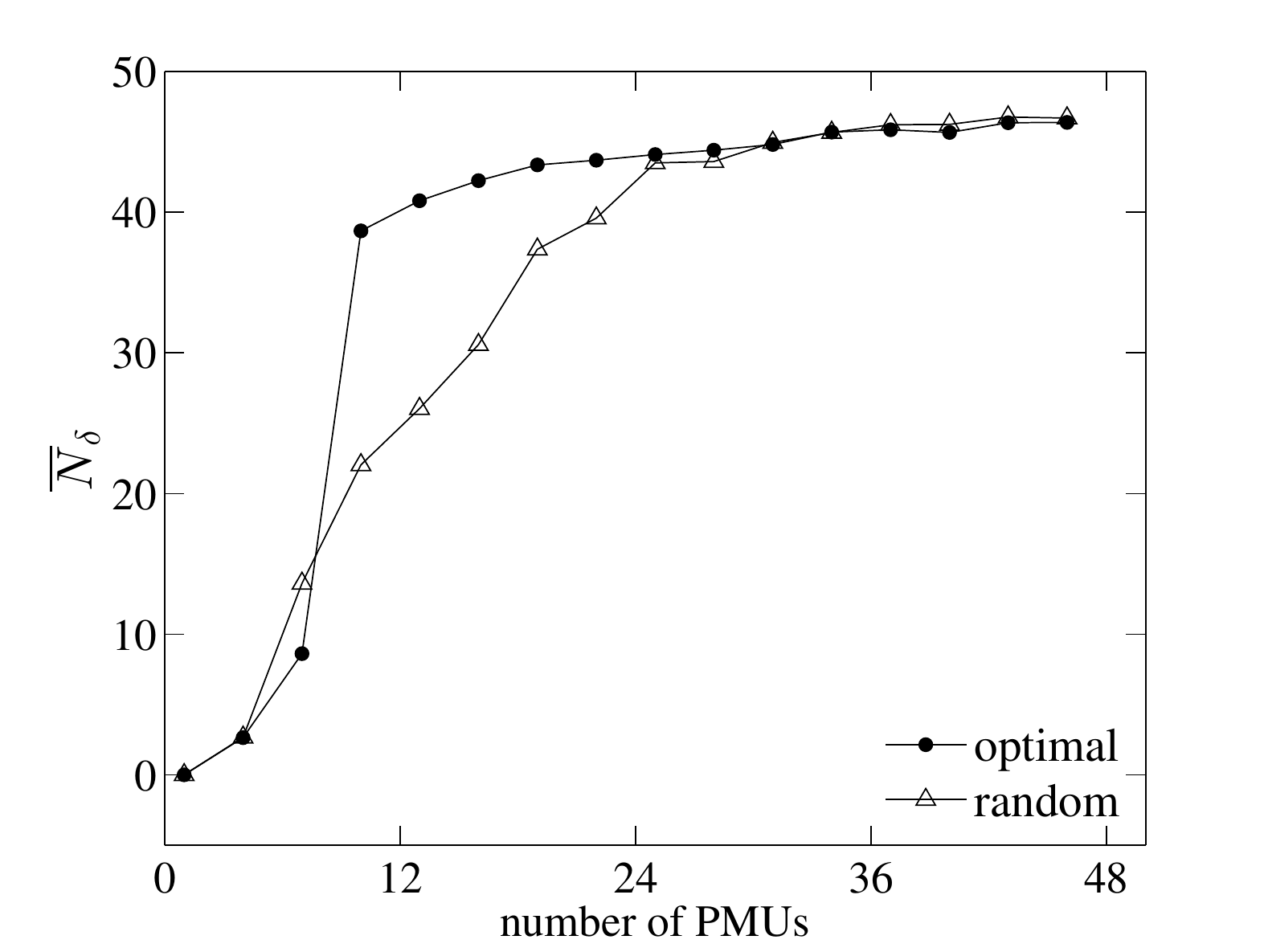}
\caption{Number of convergent angles under optimal and random PMU placements for NPCC 48-machine system.}
\label{num}
\end{figure}

The average value of the rotor angle $\bar{e}_\delta$ and the average value of the number of convergent angles $\overline{N}_\delta$ for performing dynamic state estimation for the 50 cases under the optimal and random PMU placements are shown in Figs. \ref{d} and \ref{num}.
The results for the rotor speed and the transient voltage along $q$ and $d$ axes are similar and thus are not presented here.

Compared with random PMU placements, the optimal PMU placements have much less rotor angle error and a significantly larger number of convergent angles.
Besides, under optimal PMU placements an obvious observability transition can be observed.
When the number of PMUs $\bar{g}=10$ (about 21\% of the total number of PMUs), the rotor angle error very abruptly decreases from a very high level ($\bar{e}_\delta=3.19\times 10^{18}$) to a very low level ($\bar{e}_\delta=0.69$). Correspondingly, when $\bar{g}=10$ $\overline{N}_\delta$ increases very abruptly from less than 10 to about 40.
However, under random PMU placement there is no such obvious observability transition.
The error of the rotor angle first abruptly decreases to a not very low level ($\bar{e}_\delta=1.04\times 10^5$) when $\bar{g}=16$ and then begins to gradually decreases. In order to get a very low level of rotor angle error it needs as many as 37 PMUs. Very different from what is observed for optimal PMU placement, under random PMU placement $\overline{N}_\delta$ increases approximately linearly before saturating at a high level, indicating a much slower growth rate compared with that under optimal PMU placement. Therefore, the results shown in Figs. 7 and 8 show that the proposed optimal PMU method leads to an obvious observability transition and ensures a good observability of the system states by only using a very small fraction of PMUs.

\subsection{Results on Robustness} \label{robu}

Here, we presents results on the robustness of the proposed optimal PMU placement method
under load fluctuations and contingencies for both WSCC 3-machine system and NPCC 48-machine system.

For WSCC 3-machine system 6 cases with load fluctuations are created.
In this system there are a total of 6 branches for each of which neither of the buses directly connects to a generator.
They can be ranked in descending order of the line flow as 5--7, 7--8, 6--9, 4--5, 4--6, and 8--9.
For each branch a three-phase fault is applied at its from bus and is cleared at near and remote end after $0.05s$ and $0.1s$.

For load fluctuation cases and contingency cases the $\textrm{\textbf{OPP}}^{fluc}_{\bar{g}}$ and $\textrm{\textbf{OPP}}^{cont}_{\bar{g}}$ for $\bar{g}=1,2$ can respectively be obtained.
For all 6 load fluctuation cases $\textrm{\textbf{OPP}}^{fluc}_{\bar{g}}$ are all the same as $\textrm{\textbf{OPP}}^{ty}_{\bar{g}}$, thus the average values of $\mathcal{R}_{\bar{g}}^{fluc}$ for all 6 cases are 1.0. For the 6 contingency cases only the location of placing one PMU for the second contingency is different from $\textrm{\textbf{OPP}}^{ty}_{1}$. The average values of $\mathcal{R}_{1}^{cont}$ and $\mathcal{R}_{2}^{cont}$ for all 6 cases are 0.83 and 1.0. For the only one different case $\textrm{\textbf{OPP}}^{ty}_{1}$ is $\{3\}$ but $\textrm{\textbf{OPP}}^{cont}_{1}$ is $\{2\}$. The logarithm of the empirical observability Gramian under the second contingency when placing PMU at generator 2 is 28.95, which is smaller than 30.45 for placing PMU at generator 3 but is still greater than 24.05 for placing PMU at generator 1.

\begin{figure}[!t]
\centering
\includegraphics[width=2.9in]{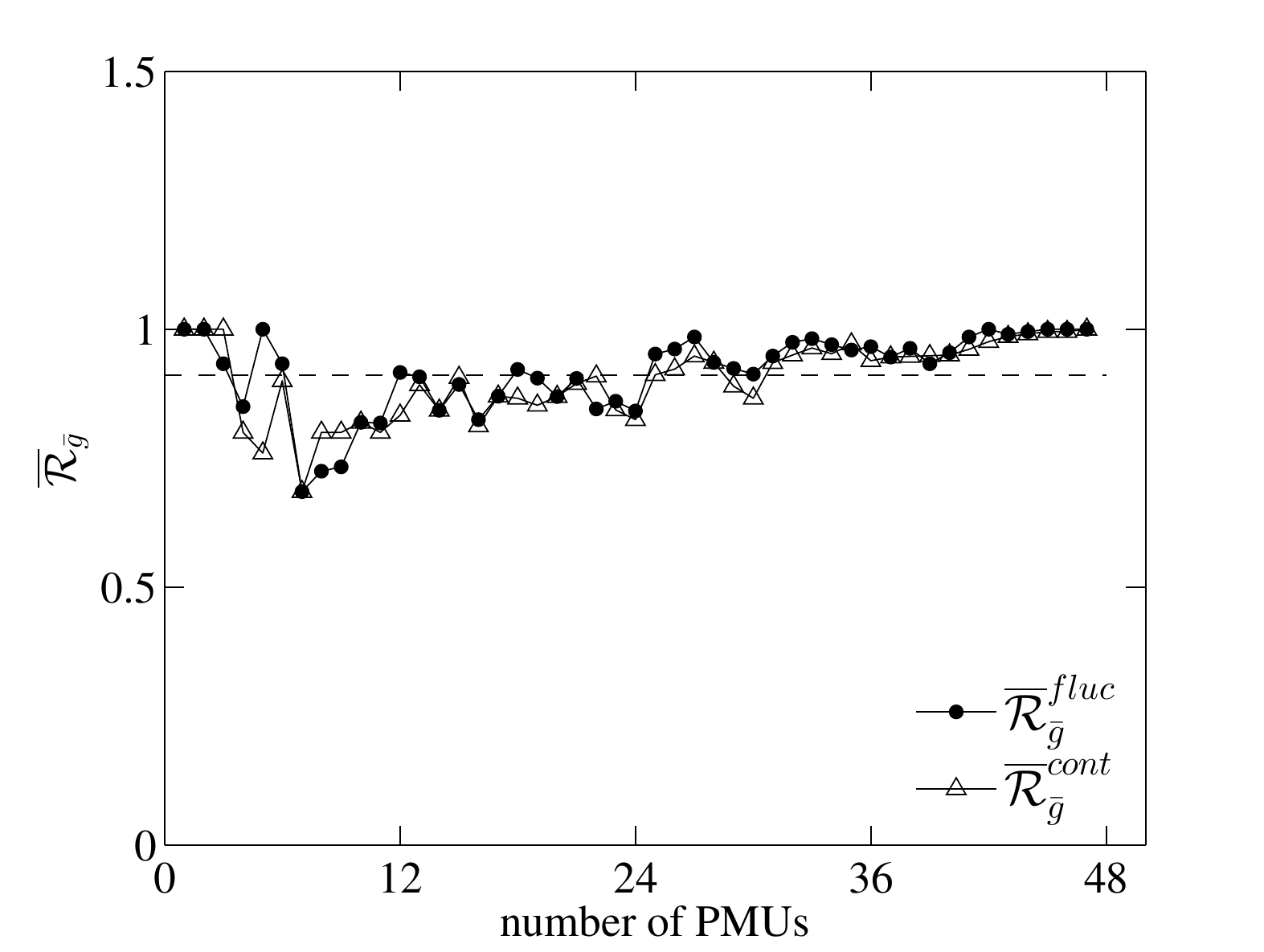}
\caption{Ratio of unchanged PMU locations under load fluctuations and contingencies for NPCC 48-machine system. The dash line indicates the average value of $\overline{\mathcal{R}}_{\bar{g}}^{fluc}$ and $\overline{\mathcal{R}}_{\bar{g}}^{cont}$ for $\bar{g}=1,\cdots,47$. }
\label{r1}
\end{figure}

Similarly, 5 cases with load fluctuations are created for NPCC 48-machine system.
Another 5 cases with contingencies are also created by applying three-phase fault at the from bus of five branches with the highest line flow.
The five branches are 132--127, 127--124, 13--12, 33--32, and 58--59 and are highlighted in Fig. \ref{npcc48}.
The average values of $\mathcal{R}_{\bar{g}}^{fluc}$ and $\mathcal{R}_{\bar{g}}^{cont}$ defined in Section \ref{robustness} for 5 cases,
which are denoted by $\overline{\mathcal{R}}_{\bar{g}}^{fluc}$ and $\overline{\mathcal{R}}_{\bar{g}}^{cont}$, are shown in Fig. \ref{r1}.
The average value of $\overline{\mathcal{R}}_{\bar{g}}^{fluc}$ and $\overline{\mathcal{R}}_{\bar{g}}^{cont}$ for $\bar{g}=1,\cdots,47$, shown by dash line in Fig. \ref{r1}, is 0.91 and is very close to 1.0.

\begin{figure}[H]
\centering
\includegraphics[width=2.9in]{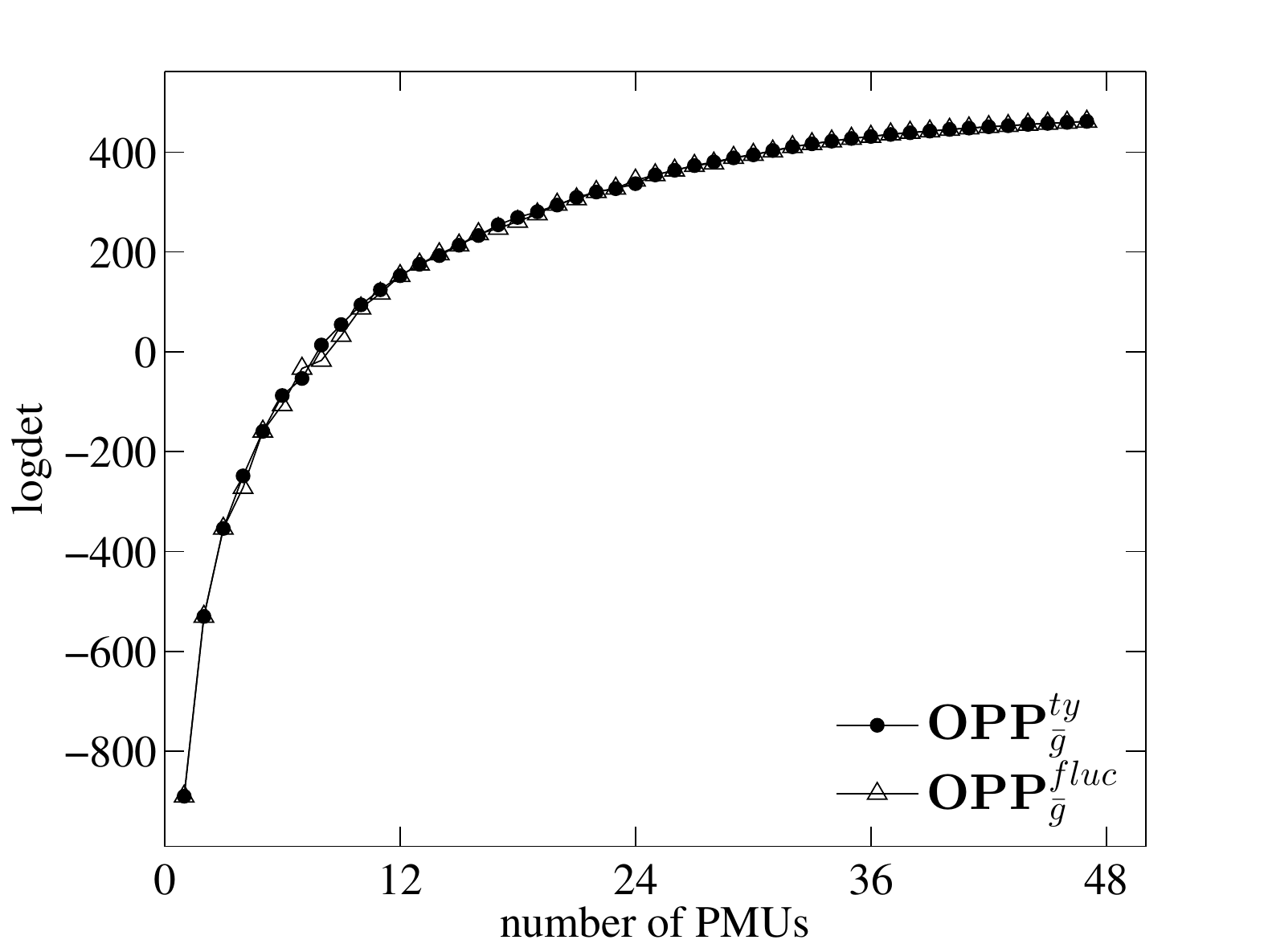}
\caption{Logarithms of determinant of the empirical observability Gramian under load fluctuations for NPCC 48-machine system.}
\label{r2}
\end{figure}

\begin{figure}[H]
\centering
\includegraphics[width=2.9in]{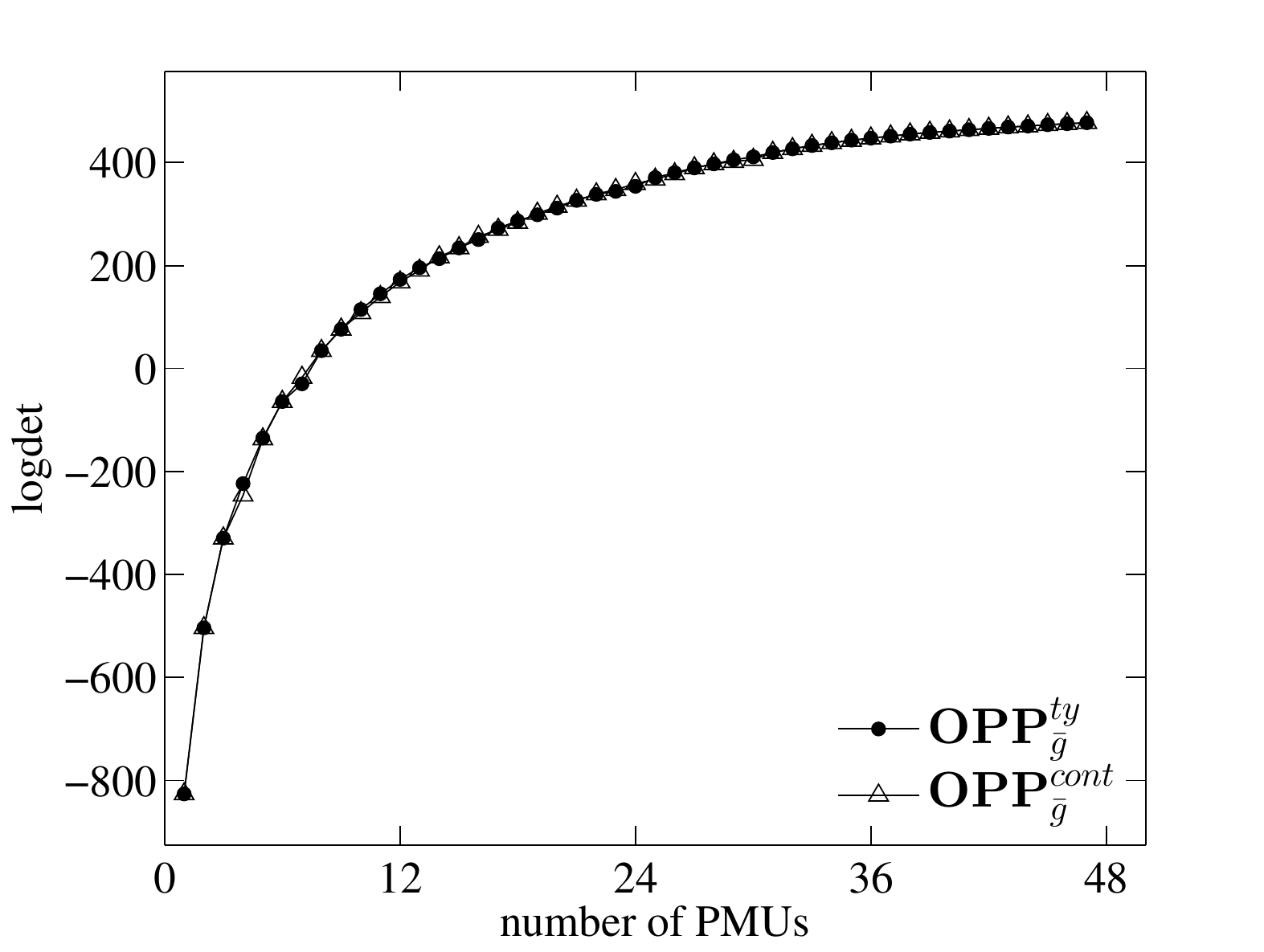}
\caption{Logarithms of determinant of the empirical observability Gramian under contingencies for NPCC 48-machine system.}
\label{r3}
\end{figure}

We also compare the logarithm of the empirical observability Gramian by using the method in Section \ref{robustness} and the results for the first load fluctuations case and the first contingency case are shown in Figs. \ref{r2} and \ref{r3}. The results for all the other cases are similar and are not given.
Note that for the first contingency case the fault is applied on one end of the line with the highest line flow, which to some extent corresponds to the most severe N-1 contingency.
It is seen that the logarithm of the empirical observability Gramian for $\textrm{\textbf{OPP}}^{ty}_{\bar{g}}$ is very close to that for $\textrm{\textbf{OPP}}^{fluc}_{\bar{g}}$ and $\textrm{\textbf{OPP}}^{cont}_{\bar{g}}$, indicating that good observability can still be guaranteed if the optimal PMU placement under typical power flow conditions is used even when there are load fluctuations or contingencies in the system.

\section{Conclusion} \label{conclusion}

In this paper the empirical observability Gramian is applied to quantify the degree of observability of a power system under a specific PMU configuration and an optimal PMU placement method for dynamic state estimation
is proposed by maximizing the determinant of the empirical observability Gramian.
It is effectively and efficiently solved by the NOMAD solver and is then tested on WSCC 3-machine 9-bus system and
NPCC 48-machine 140-bus system by performing dynamic state estimation with square-root unscented Kalman filter.

The results show that the obtained optimal PMU placements can guarantee smaller estimation errors and
larger number of convergent states compared with random PMU placements.
Under optimal PMU placements an obvious observability transition can be observed.
Although the optimal PMU placements is obtained for the system under typical power flow conditions the obtained optimal placements are very robust to load fluctuations and contingencies and can still guarantee good observability even under small or big disturbances.


%

\section*{Acknowledgment}

We gratefully thank anonymous reviewers for their insightful advice
that helped greatly improve our paper.

\begin{IEEEbiography} [{\includegraphics[width=1in,height=1.25in,clip,keepaspectratio]{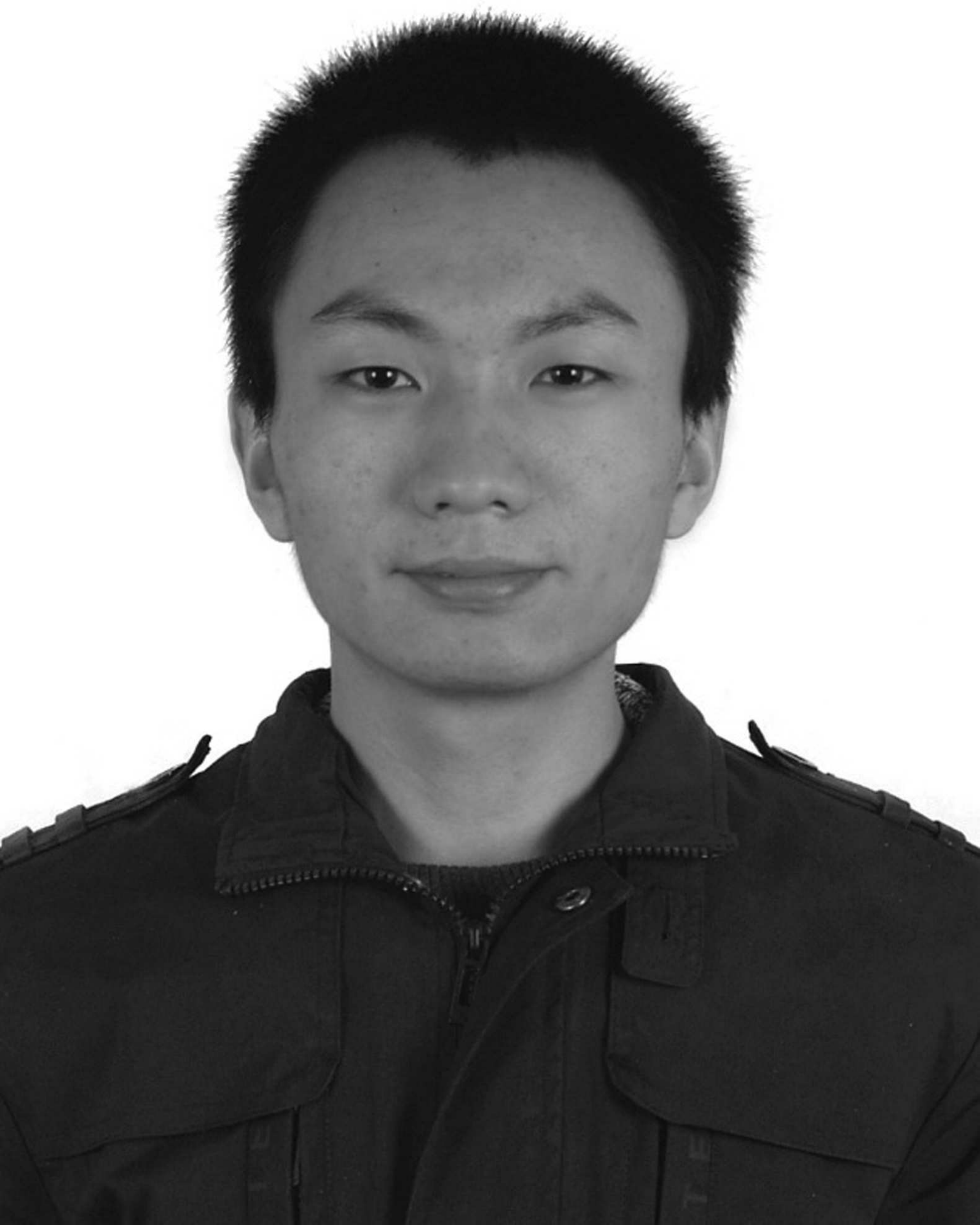}\vfill}]
{Junjian Qi} (S'12--M'13)
received the B.E. and Ph.D. degree both in Electrical Engineering from Shandong University in 2008 and Tsinghua University in 2013.

He visited Prof. Ian Dobson's group at Iowa State University in Feb.--Aug. 2012 and is currently a Research Associate at Department of Electrical Engineering and Computer Science, University of Tennessee, Knoxville, TN, USA. His research interests include blackouts, cascading failure, state estimation, and synchrophasors.
\end{IEEEbiography}

\begin{IEEEbiography} [{\includegraphics[width=1in,height=1.25in,clip,keepaspectratio]{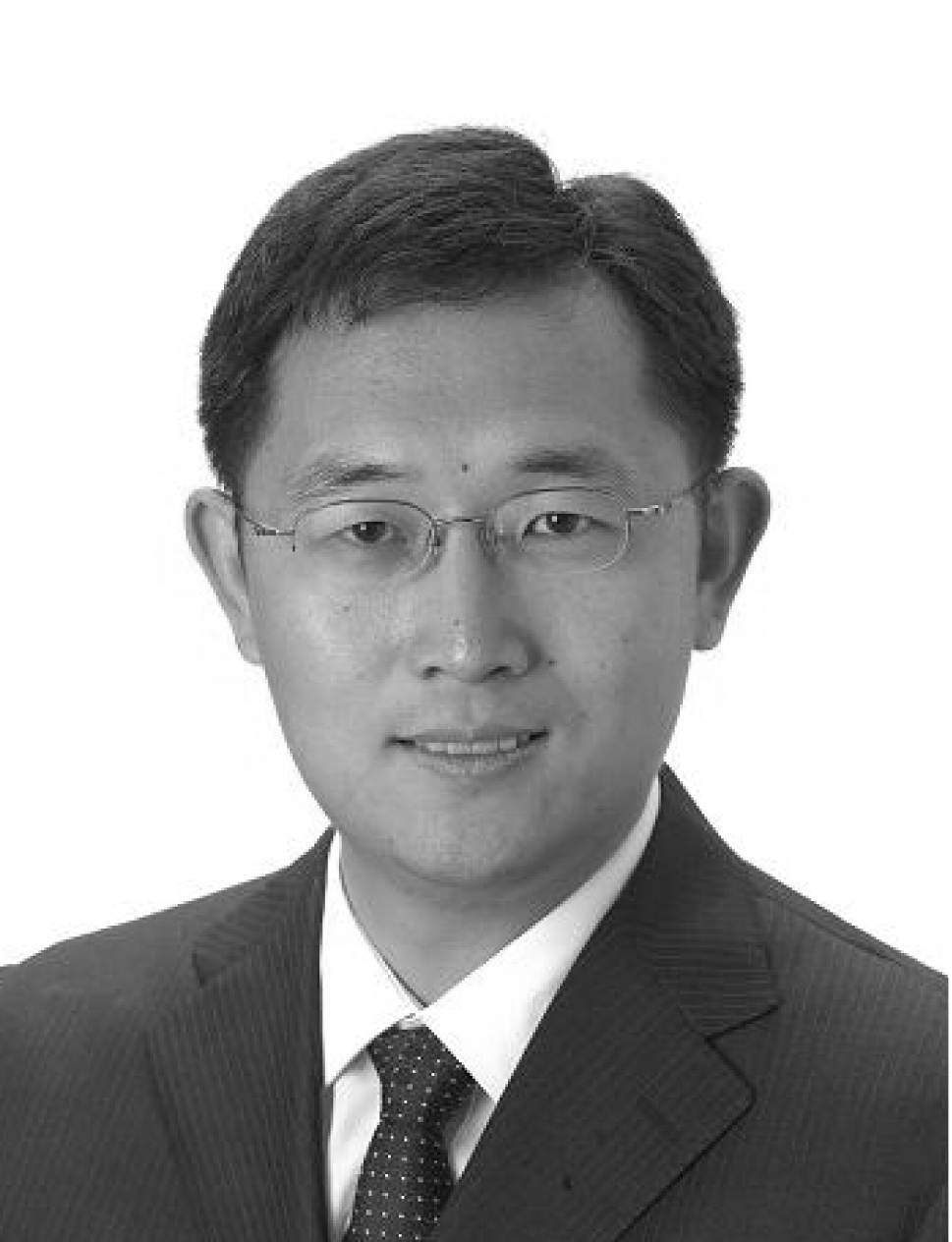}\vfill}]
{Kai Sun} (M'06--SM'13)
received the B.S. degree in automation in 1999 and the Ph.D. degree in control
science and engineering in 2004  both  from Tsinghua University, Beijing, China.

He was a Postdoctoral Research Associate at Arizona State University, Tempe, from 2005 to 2007,
and was a Project Manager in grid operations and planning areas at EPRI, Palo Alto, CA
from 2007 to 2012.
He is currently an Assistant Professor at the Department of Electrical Engineering and Computer Science, University of Tennessee, Knoxville, TN, USA. He serves as an editor in IEEE Transactions on Smart Grid.
\end{IEEEbiography}

\begin{IEEEbiography} [{\includegraphics[width=1in,height=1.25in,clip,keepaspectratio]{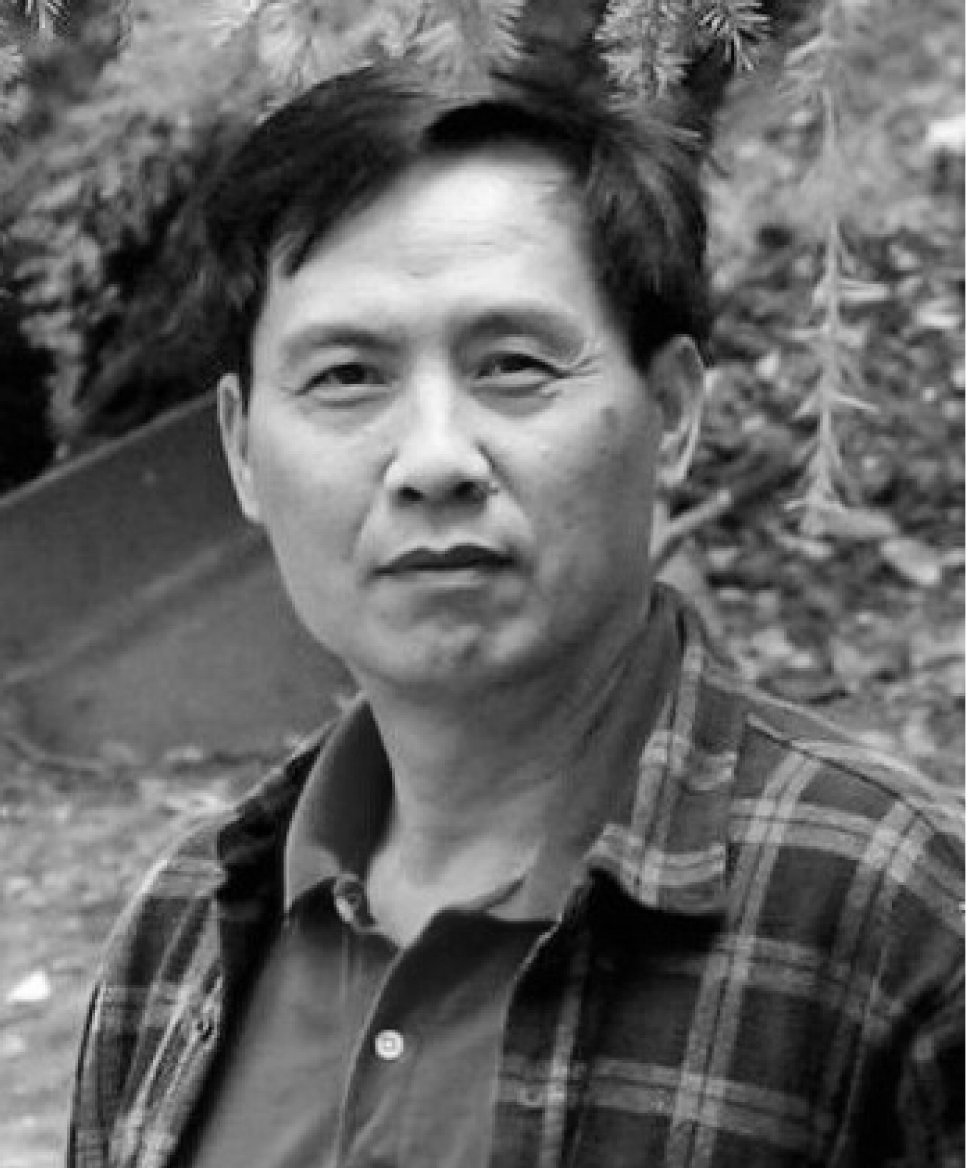}\vfill}]
{Wei Kang} (M'91--F'08)
received the B.S. and M.S. degrees from Nankai University, China, both in mathematics, in 1982 and 1985, respectively, and the Ph.D. degree in mathematics from the University of California, Davis, in 1991.

He is currently a Professor of applied mathematics at the U.S. Naval Postgraduate School, Monterey, CA.
He was a visiting Assistant Professor of systems science and mathematics at Washington University, St. Louis, MO (1991-1994).
He served as the Director of Business and International Collaborations at the American Institute of Mathematics (2008-2011).
His research interest includes computational optimal control, nonlinear filtering, cooperative control of autonomous vehicles, industry applications of control theory, nonlinear $H_\infty$ control, and bifurcations
and normal forms. His early research includes topics on Lie groups, Lie algebras, and differential geometry.

Dr. Kang is a fellow of IEEE. He was a plenary speaker in several international conferences of SIAM and IFAC. He served as an associate editor in several journals, including IEEE TAC and Automatica.
\end{IEEEbiography}


\begin{thebibliography}{11}

\bibitem{se1} \label{se1}
F. C. Schweppe and J. Wildes, ``Power system static-state estimation, Part \uppercase\expandafter{\romannumeral 1\relax}: exact model," \emph{IEEE Trans. Power App. Syst.}, vol. PAS-89, no. 1, pp. 120--125, Jan. 1970.

\bibitem{se2} \label{se2}
A. Abur and A. G\'{o}mez Exp\'{o}sito, \emph{Power System State Estimation: Theory and Implementation}. CRC Press, 2004.

\bibitem{se3} \label{se3}
A. Monticelli, ``Electric power system state estimation," \emph{Proc. IEEE}, vol. 88, no. 2, pp. 262--282, Feb. 2000.

\bibitem{se4} \label{se4}
M. R. Irving, ``Robust state estimation using mixed integer programming," \emph{IEEE Trans. Power Syst.}, vol. 23, no. 3, pp. 1519--1520, Aug. 2008.

\bibitem{se5} \label{se5}
G. He, S. Dong, J. Qi, and Y. Wang, ``Robust state estimator based on maximum normal measurement rate," \emph{IEEE Trans. Power Syst.}, vol. 26, no. 4, pp. 2058--2065, Nov. 2011.

\bibitem{se6} \label{se6}
J. Qi, G. He, S. Mei, and F. Liu, ``Power system set membership state estimation," in \emph{Proc. IEEE Power and Energy Society General Meeting}, pp. 1--7, San Diego, CA USA, Jul. 2012.

\bibitem{sun1} \label{sun1}
K. Sun, K. Hur, and P. Zhang, ``A new unified scheme for controlled power system separation using synchronized phasor measurements," IEEE Transactions on Power Systems, vol. 26, no. 3, pp. 1544--1554, Aug. 2011

\bibitem{sun2} \label{sun2}
K. Sun, Q. Zhou, and Y. Liu, ``A Phase Locked Loop-based Approach to Real-time Modal Analysis on Synchrophasor Measurements," IEEE Transactions on Smart Grid, vol. 5, no. 1, pp. 260--269, Jan. 2014

\bibitem{kailath} \label{kailath}
T. Kailath, \emph{Linear Systems}, Prentice-Hall: Englewood Cliffs, NJ, 1980.

\bibitem{diop} \label{diop}
S. Diop and M. Filess, ``On nonlinear observability," \emph{Proceedings of ECC'91}, Herm\`{e}s, Paris, vol. 1, pp. 152--157, 1991.

\bibitem{diop1} \label{diop1}
S. Diop and M. Filess, ``Nonlinear observability, identifiability, and persistent trajectories," \emph{Proceedings of the 30th IEEE Conference on Decision and Control}, vol. 1, pp. 714--719, 1991.

\bibitem{lall} \label{lall}
S. Lall, J. E. Marsden, and S. Glava\v ski, ``Empirical model reduction of controlled nonlinear systems," \emph{14th IFAC World Congress, Beijing China}, pp. 473--478, 1999.

\bibitem{lall1} \label{lall1}
S. Lall, J. E. Marsden, and S. Glava\v ski. ``A subspace approach to balanced truncation for model reduction of nonlinear control systems," \emph{International Journal of Robust and Nonlinear Control}, vol. 12, pp. 519--535, 2002.

\bibitem{Gramianref} \label{Gramianref}
J. Hahn and T. F. Edgar, ``Balancing approach to minimal realization and model reduction of stable nonlinear systems," \emph{Industrial and Engineering Chemistry Research}, vol. 41, no. 9, pp. 2204–-2212, 2002.

\bibitem{Krener} \label{Krener}
A. J. Krener and K. Ide, ``Measures of unobservability," \emph{IEEE Conference on Decision and Control}, Shanghai, China, pp. 6401--6406, Dec. 2009.

\bibitem{kang1} \label{kang1}
W. Kang and L. Xu, ``Computational analysis of control systems using dynamic optimization," \emph{arXiv:0906.0215v2}, 2009.

\bibitem{sun} \label{sun}
K. Sun and W. Kang, ``Observability and estimation methods using sychrophasors," in \emph{Proc. IFAC World Congr.}, pp. 963--968, 2014.

\bibitem{kang} \label{kang}
W. Kang and L. Xu, ``Optimal placement of mobile sensors for data assimilations," \emph{Tellus A}, 64, 17133, 2012.

\bibitem{singh} \label{singh}
A. K. Singh and J. Hahn, ``Determining optimal sensor locations for state and parameter estimation for stable nonlinear systems," \emph{Ind. Eng. Chem. Res.}, vol. 44, no. 15, pp. 5645--5659, 2005.

\bibitem{singh1} \label{singh1}
A. K. Singh and J. Hahn, ``Sensor location for stable nonlinear dynamic systems: multiple sensor case," \emph{Ind. Eng. Chem. Res.}, vol. 45, no. 10, pp. 3615--3623, 2006.

\bibitem{det} \label{det}
M. Serpas, G. Hackebeil, C. Laird, and J. Hahn, ``Sensor location for nonlinear dynamic systems via observability analysis and MAX-DET optimization," \emph{Computers \& Chemical Engineering}, vol. 48, pp. 105--112, 2012.

\bibitem{ob1} \label{ob1}
T. Baldwin, L. Mili, J. Boisen, M.B., and R. Adapa, ``Power system observability with minimal phasor measurement placement," \emph{IEEE Trans. Power Syst.}, vol. 8, no. 2, pp. 707--715, May 1993.

\bibitem{info} \label{info}
Q. Li, T. Cui, Y. Weng, R. Negi, F. Franchetti, and M. Ili\'c, ``An information-theoretic approach to PMU placement in electric power systems," \emph{IEEE Trans. Smart Grid}, vol. 4, pp. 446--456, Mar. 2013.

\bibitem{ob2} \label{ob2}
B. Xu and A. Abur, ``Observability analysis and measurement placement for systems with PMUs," \emph{in IEEE Power Syst. Conf. Expo.}, vol. 2, pp. 943--946, Oct. 2004.

\bibitem{ob3} \label{ob3}
B. Gou, ``Optimal placement of PMUs by integer linear programming," \emph{IEEE Trans. Power Syst.}, vol. 23, no. 3, pp. 1525--1526, Aug. 2008.

\bibitem{ob4} \label{ob4}
S. Chakrabarti and E. Kyriakides, ``Optimal placement of phasor measurement units for power system observability," \emph{IEEE Trans. Power Syst.}, vol. 23, no. 3, pp. 1433--1440, Aug. 2008.

\bibitem{ob5} \label{ob5}
B. Milosevic and M. Begovic, ``Nondominated sorting genetic algorithm for optimal phasor measurement placement," \emph{IEEE Trans. Power Syst.}, vol. 18, no. 1, pp. 69--75, Feb. 2003.

\bibitem{ob6} \label{ob6}
F. Aminifar, C. Lucas, A. Khodaei, and M. Fotuhi-Firuzabad, ``Optimal placement of phasor measurement units using immunity genetic algorithm," \emph{IEEE Trans. Power Del.}, vol. 24, no. 3, pp. 1014--1020, Jul. 2009.

\bibitem{ob7} \label{ob7}
S. Chakrabarti, G. K. Venayagamoorthy and E. Kyriakides, ``PMU placement for power system observability using binary particle swarm optimization," \emph{Australasian Universities Power Engineering Conference}, Brisbane, Australia, pp. 1--5, Dec. 2008.

\bibitem{ob8} \label{ob8}
A. Almutairi and J. Milanovi, ``Comparison of Different Methods for Optimal Placement of PMUs," \emph{Proceedings of 2009 IEEE Bucharest Power Tech Conference}, Bucharest, Romania, pp. 1--6, Jun. 2009.

\bibitem{ob9} \label{ob9}
P. Korba, M. Larsson and C. Rehtanz, ``Detection of Oscillations in Power Systems using Kalman Filtering Techniques," \emph{Proceedings of the IEEE Conference on Control Applications}, Istanbul, Turkey, pp. 183--188, Jun. 2003.

\bibitem{kamwa} \label{kamwa}
I. Kamwa and R. Grondin, ``PMU configuration for system dynamic performance measurement in large multiarea power systems," \emph{IEEE Trans. Power Syst.}, vol. 17, no. 2, pp. 385--394, May 2002.

\bibitem{zhang1} \label{zhang1}
J. Zhang, G. Welch, and G. Bishop, ``Observability and estimation uncertainty analysis for PMU placement alternatives," \emph{Proceedings of North American Power Symposium}, Arlington, TX, USA, pp. 1--8, Sep. 2010.

\bibitem{zhang2} \label{zhang2}
J. Zhang, G. Welch, G. Bishop, and Z. Huang, ``Optimal PMU placement evaluation for power system dynamic state estimations," \emph{Proceedings of IEEE PES Conference on Innovative Smart Grid Technologies Europe}, G\"{o}teborg, Sweden, pp. 1--7, Oct. 2010.

\bibitem{huang} \label{huang}
Y. Sun, P. Du, Z. Huang, K. Kalsi, R. Diao, K. K. Anderson, Y. Li, and B. Lee, ``PMU placement for dynamic state tracking of power systems," \emph{North American Power Symposium}, Boston MA, pp. 1--7, Aug. 2011.

\bibitem{muller} \label{muller}
P. C. M\'{u}ller and H. I. Weber, ``Analysis and optimization of certain quantities of controllability and observability for linear dynamical systems," \emph{Automatica}, vol. 8, no. 3, pp. 237--246, May 1972.

\bibitem{zhou} \label{zhou}
N. Zhou, D. Meng, and S. Lu, ``Estimation of the dynamic states of synchronous machines using an extended particle filter," \emph{IEEE Trans. Power Syst.}, vol. 28, no. 4, pp. 4152--4161, Nov. 2013.

\bibitem{kunder} \label{kunder}
P. Kunder, \emph{Power System Stability and Control}, New York, NY, USA: McGraw-Hill, 1994.

\bibitem{emgr} \label{emgr}
C. Himpe and M. Ohlberger, ``A unified software framework for empirical gramians," \emph{Journal of Mathematics}, vol. 2013, pp. 1--6, 2013.

\bibitem{boyd} \label{boyd}
L. Vandenberghe, S. Boyd, and S. Wu, ``Determinant maximization with linear matrix inequality constraints," \emph{SIAM J. Matrix Anal. Appl.}, vol. 19, no. 2, pp. 499--533, 1998.

\bibitem{wang_det} \label{wang_det}
C. Wang, D. Sun, and K. Toh, ``Solving log-determinant optimization problems by a Newton-CG primal proximal point algorithm," \emph{SIAM J. Optim.}, vol. 20, no. 6, pp. 2994--3013, 2010).

\bibitem{nomad} \label{nomad}
S. Le Digabel, ``Algorithm 909: NOMAD: Nonlinear optimization with the MADS algorithm," \emph{ACM Trans. Mathematical Software}, vol. 37, no. 4, Feb. 2011.

\bibitem{opti} \label{opti}
J. Currie and D. I. Wilson, ``OPTI: lowering the barrier between open source optimizers and the industrial MATLAB user," \emph{Foundations of Computer-Aided Process Operations}, Georgia, USA, pp. 8--11, 2012.

\bibitem{mads} \label{mads}
C. Audet and J. E. Dennis Jr., ``Mesh adaptive direct search algorithms for constrained optimization,"
\emph{SIAM J. Optim.}, vol. 17, no. 1, pp. 188--217, 2006.

\bibitem{clarke} \label{clarke}
F. H. Clarke, \emph{Optimization and Nonsmooth Analysis}. New York, NY, USA: Wiley, 1983.

\bibitem{audet} \label{audet}
C. Audet, V. B\'{c}hard, and S. Le Digabel, ``Nonsmooth optimization through mesh adaptive direct search and variable neighborhood search," \emph{J. Global Optim.}, vol. 41, no. 2, pp. 299--318, 2008.

\bibitem{vns} \label{vns}
N. Mladenovi\'{c} and P. Hansen, ``Variable neighborhood search," \emph{Computer \& Operations Research}, vol. 24, no. 11, pp. 1097--1100, Nov. 1997.

\bibitem{sr_ukf} \label{sr_ukf}
R. Merwe and E. Wan, ``The square-root unscented Kalman filter for state and parameter-estimation," \emph{Proc. IEEE Int. Conf. Acoustics, Speech, and Signal Processing (ICASSP)}, vol. 6, pp. 3461--3464, 2001.

\bibitem{CarrerasHICSS13} \label{OPA}
B. A. Carreras, D. E. Newman, I. Dobson, and N. S. Degala, ``Validating OPA with WECC data," \emph{46th Hawaii Intl. Conference on System Sciences}, Maui, HI, pp. 2197--2204, Jan. 2013.

\bibitem{bp13} \label{bp13}
J. Qi, I. Dobson, and S. Mei, ``Towards estimating the statistics of simulated cascades of outages with branching processes,'' \emph{IEEE Trans. Power Syst.}, vol. 28, no. 3, pp. 3410--3419, Aug. 2013.

\bibitem{greedy} \label{greedy}
A. Conejo, T. G\'{o}mez, and J. I. de la Fuente, ``Pilot-bus selection for secondary voltage control," \emph{Eur. Trans. Elec. Eng.}, vol. 5, no. 5, pp. 359--366, 1993.

\bibitem{algo} \label{algo}
P. Bonami, L. T. Biegler, A. R. Conn, et al, ``An algorithmic framework for convex mixed integer nonlinear
programs.Discrete Optimization," \emph{Discrete Optimization}, vol. 5, no. 2, pp. 186--204, May 2008.

\bibitem{anderson} \label{anderson}
P. M. Anderson and A. A. Fouad, \emph{Power System Control and Stability}, Iowa State University Press, Ames Iowa, 1977.

\bibitem{pst} \label{pst}
J. Chow and G. Rogers, User manual for power system toolbox, version 3.0, 1991--2008.

\bibitem{huang1} \label{huang1}
Z. Huang, P. Du, D. Kosterev, and B. Yang, ``Application of extended Kalman filter techniques for dynamic model parameter calibration," \emph{IEEE  Power  Engineering  Society  General  Meeting}, Calgary, Canada, pp. 1--8, Jul. 2009.


\end{thebibliography}
\end{document}